\documentclass[lean]{Draft}

\usepackage{newunicodechar}
\usepackage{xspace}
\usepackage{amsmath, amssymb}
\usepackage{colonequals}
\usepackage[colorlinks=true]{hyperref}
\usepackage{todonotes}

\newunicodechar{⋃}{\cup}

\newcommand*{\R}{\mathbb{R}}

\newcommand*{\mathlib}{\textsc{Mathlib}\xspace}
\newcommand*{\Lean}{\textsc{Lean}\xspace}

\usepackage{tikz}
\usetikzlibrary{arrows, arrows.meta}
\definecolor{ghpurple}{rgb}{0.4, 0.22, 0.73}

\renewcommand{\lean}[1]{\lstinline| #1|}

\makeatletter
\newcommand\leanlink{\begingroup\catcode`\#=12\relax\@leanlink}
\newcommand\@leanlink[2]{\endgroup
\href{#1}
{\texttt{\detokenize{#2}}}}

\newcommand{\mllink}[3]{%
\leanlink{https://github.com/leanprover-community/mathlib4/blob/82d756ea359cfa0f2c10bc00b8a997822c3b32ea/Mathlib/#1.lean\##2}
{#3}}

\newcommand{\kelink}[3]{%
\leanlink{https://github.com/RemyDegenne/kolmogorov_extension4/blob/master/KolmogorovExtension4/#1.lean\##2}
{#3}}

\newcommand{\bmlink}[3]{%
\leanlink{https://github.com/RemyDegenne/brownian-motion/blob/6a7a3141c6bc217484f635f3ab861316a37c0362/BrownianMotion/#1.lean\##2}
{#3}}

\title[Formalization of Brownian motion in Lean]{Formalization of Brownian motion in Lean}
\author[R. Degenne, D. Ledvinka, E. Marion, P. Pfaffelhuber]{
  Rémy Degenne, David Ledvinka, Etienne Marion, Peter Pfaffelhuber}

\authorinfo[R. Degenne]{Univ. Lille, Inria, CNRS, Centrale Lille, UMR 9189-CRIStAL, F-59000 Lille, France}{remy.degenne@inria.fr}
\authorinfo[D. Ledvinka]{University of Toronto, Canada}{david.ledvinka@mail.utoronto.ca}
\authorinfo[E. Marion]{ENS de Lyon, 46, Allée d’Italie, 69007 Lyon, France}{etienne.marion@ens-lyon.fr}
\authorinfo[P. Pfaffelhuber]{University of Freiburg, Germany}{p.p@stochastik.uni-freiburg.de}

\msc{TODO}
\keywords{Formalization, Mathlib, Lean, Probability, Brownian motion}

\VOLUME{TODO}
\YEAR{TODO}
\NUMBER{TODO}
\firstpage{1}
\DOI{TODO}
\receiveddate{TODO}
\finaldate{TODO}
\accepteddate{TODO}

\begin{abstract}
Brownian motion is a building block in modern probability theory. In this paper, we describe a formalization of Brownian motion using the \Lean theorem prover. We build on the existing measure-theoretic foundations in Lean's mathematical library, Mathlib, and we develop several key components needed for the construction of Brownian motion, including the Carathéodory and Kolmogorov extension theorems, Gaussian measures in Banach spaces, and the Kolmogorov-Chentsov theorem for path continuity.
\end{abstract}

\begin{document}




\section{Introduction}

\subsection{Mathematical background}
\sloppy Brownian motion is arguably one of the most important stochastic processes (e.g.\ \cite{karatzas1991brownian, morters2010brownian}), and is used as a modeling tool across all sciences (physics: e.g.\ \cite{einstein1906theorie, bian2016111}; biology: e.g.\ \cite{erban2014molecular}; finance: e.g.\ \cite{davis2006louis}). Mathematically, Brownian motion led to the Wiener measure \cite{wiener1923differential}, which is the first instance of a probability measure on a function space, but also to developments such as Stochastic (Partial) Differential equations (see e.g.\ \cite{hairer2009introduction}).

The goal of the present paper is to describe a formalization of Brownian motion using \Lean \cite{lean, moura2021lean}, building on its mathematical library \mathlib \cite{mathlib}. For this paper we use the following definition of Brownian motion (there are many equivalent characterizations, e.g.\ Definition 2.12 with Proposition 2.3 of~\cite{LeGall2016}).

\begin{definition}[\bmlink{Gaussian/BrownianMotion}{L490-L493}{ProbabilityTheory.IsBrownian}]\label{def:IsBrownian_preview}
A stochastic process $X$ on $\mathbb{R}_+$ with values in $\mathbb{R}$ is said to be a \emph{Brownian motion} if $X$ has continuous paths almost surely and for all $n \in \mathbb{N}$ and $t_1, \ldots, t_n \in \mathbb{R}_+$, the vector $(X_{t_1}, \ldots, X_{t_n})$ has a multivariate normal distribution with mean $0$ and covariance matrix given by $\mathbf{cov}(X_{t_i}, X_{t_j}) = t_i \wedge t_j$ (the minimum of $t_i$ and $t_j$).
\end{definition}
Note that it is a non-trivial fact that a process exists which satisfies the above definition, and most of the formalization is dedicated to constructing such a process (or developing the necessary prerequisites)
and proving it satisfies these properties. Here are some other interesting facts about Brownian motion which we have formalized:
\begin{itemize}
\item $X_t$ has almost surely locally $\beta$-Hölder-continuous paths for any coefficient $\beta < \tfrac 12$.
\item For any $t_0 \in \R_+$, the process $X_{t_0 + t} - X_{t_0}$ is a Brownian motion independent of ${(X_t)}_{t \leq t_0}$. In particular $X_t$ satisfies the weak Markov property.
\item For $c > 0$, the process $X_{ct}/\sqrt{c}$ is a Brownian motion.
\item $tX_{1/t}$ is a Brownian motion.
\item $X_t/t$ tends to $0$ almost surely as $t$ tends to infinity.
\end{itemize}
And here are some which we leave for future work:
\begin{itemize}
\item $X_t$ is almost surely nowhere differentiable.
\item For every almost surely finite stopping time $T$, the process $X_{T + t} - X_{T}$ is a Brownian motion independent of $\mathcal{F}_T$ (the $\sigma$-algebra of events determined by time $T$). In particular $X_t$ satisfies the strong Markov property.
\item Brownian motion is a martingale, i.e.\ $\mathbb E[X_{t+s} | (X_r)_{r\leq t}] = X_t$ for $s,t \geq 0$; moreover it is the only stochastic process with continuous paths such that $(X_t)_{t\geq 0}$ and $(X_t^2 - t)_{t\geq 0}$ are both martingales.
\end{itemize}

At the start of the project, \mathlib had a solid foundation in measure theory but lacked other tools we needed like multivariate Gaussian distributions. Starting from what was available, we arrived at a fully formalized construction of Brownian motion on $\mathbb{R}_+$ by following these steps:
\begin{itemize}
\item {\em The Carathéodory and Kolmogorov Extension Theorems.} Constructing a probability measure on an uncountable infinite product (i.e.\ the distribution of a stochastic process) is not straightforward, since the corresponding product of $\sigma$-algebras is too large.
Rather, one uses the $\sigma$-algebra defined by all finite (or countably infinite) projections, and constructs a probability measure by extension from all finite dimensional possibilities, which is the content of the Kolmogorov Extension Theorem. This result is in fact an application of a more general result, the Carathéodory Extension Theorem. In this step, we need that the family of finite dimensional distributions has a consistency property (i.e.\ it forms a projective family).
\item {\em Gaussian Measures and characteristic functions.} Defining a one-dimensional normal distribution is straight-forward using its density with respect to Lebesgue measure, and concrete calculations are often possible using the density. In the multi-dimensional case, which we have to work with extensively, it is often much more convenient to use characteristic functions, which are known to characterize probability measures uniquely.
In particular, they can be used to show that linear maps (e.g.\ projections) of Gaussian measures are Gaussian. This approach is fundamental to show that the finite dimensional distributions of the Brownian motion indeed form a projective family.
\item {\em The Kolmogorov-Chentsov Theorem.} We are aiming for a stochastic process with continuous paths. The classical Kolmogorov Chentsov theorem gives a criterion in terms of some moment bounds if the set of times is a subset of $\mathbb R^d$; see e.g.\ \cite{kallenberg2021}. Here, we use a modern version for existence of a modification of a stochastic process with continuous paths, based on a more general set of times \cite{kratschmer2023kolmogorov}.
\item {\em Construction of a Brownian motion and Wiener measure on $\R_+$.} Putting all of the above together gives a continuous Gaussian stochastic process with the correct distribution, i.e.\ Brownian motion.
We show important properties of Brownian motion, such as the independent increments property.
From this process, we can also define a probability measure on the (Polish) space of continuous functions $\mathbb R_+ \to \mathbb R$. This is usually called the Wiener measure.
\end{itemize}

It is important to realize that probability theory comes with a duality between random variables (or stochastic processes) and their distributions (or laws), which are probability measures. In the case of stochastic processes, the latter are defined on the (often uncountable infinite) product on the state space. While some properties of a stochastic process are only dependent on their distribution (e.g.\ being a Gaussian process), other properties depend on more refined properties (e.g.\ having continuous paths).
Our formalization led to the addition to \mathlib of new definitions and results that allow to manipulate random variables directly, without having to perform operations on their distributions all the time; see Section~\ref{S:RV-formalism} for more details.

\subsection{Related work}
Stochastic processes have been formalized in several theorem provers.
Within Isabelle/HOL, the Kolmogorov extension theorem has been previously formalized \cite{Immler2012}. This formalization only works on Polish spaces (rather than on spaces where every finite measure is inner regular with respect to compact sets, see below), and only in the case where the state spaces for all times are identical. As for the Kolmogorov-Chentsov Theorem (showing continuity of paths), there is an Isabelle/HOL formalization as well \cite{Kolmogorov_Chentsov-AFP}.
It uses the index set  $\R_+$ (rather than a more general type) with the usual dyadics proof, as e.g. outlined in \cite{kallenberg2021}. Last, an implementation of
Brownian motion for Isabelle/HOL is described at \cite{laursen2024brownian} \footnote{the code, which is on \href{https://github.com/cplaursen/Brownian_Motion}{github}, seems to be a work in progress}.

For the Lean theorem prover, \cite{ying2023formalization} formalized theorems about discrete time martingales (see e.g.\ \cite{kallenberg2021}) and was the first implementation of stochastic processes.
Our work shares little with that formalization, apart from the general setting of stochastic processes. We however depend extensively on the measure theory content of \cite{mathlib}.

\subsection{Project organization}

This project started as a traditional collaboration between researchers, in which two of the authors (R.D. and P.P.) formalized the Kolmogorov extension theorem.
In a second phase however, we decided to transform the project into a collaborative effort, open to anyone interested in contributing.
We wrote first a succinct description of the proof we intended to follow for the Brownian motion definition, before expanding it into a \href{https://remydegenne.github.io/brownian-motion/blueprint/}{blueprint}, a collaboration tool introduced for the Sphere Eversion project \cite{sphere_eversion} (used here through the LeanProject github template \cite{Monticone_LeanProject_2025}).
This is a website that accompanies the github repository on which the code is hosted, and that contains pages with detailed descriptions of the sequence of lemmas and definitions needed to prove the main theorems.
The dependency of each lemma on the previous ones is made explicit, and the blueprint tool automatically generates a color-coded graph of dependencies.
The graph is useful to visualize the structure of the proof, and to identify which lemmas are proved, which are ready to be proved, and which are still missing prerequisites.

Before making the project public, we thus spent time detailing every step of the proofs, cutting them into small lemmas, and identifying all the missing pieces between what was available in \mathlib and what we would need for the project.
We started by implementing in \Lean several key definitions and statements (without proofs): there are always several ways to turn a mathematical concept into a \Lean definition, and in order to collaborate on a project it is important to decide early on a unique way to talk about any concept to ensure compatibility of the code produced, even if the definition chosen can then be refined.
Different parts of the blueprint had different levels of details: we wrote a precise description of the technical Kolmogorov-Chentsov proof, but omitted many proofs for more standard results on Gaussian distributions.

The project was announced on the \href{https://leanprover.zulipchat.com/}{Lean Zulip website}, and anyone interested was invited to contribute. D.L. and E.M. joined the project at that time and then contributed significantly to the formalization. The coordination of the project was done through discussions in a Zulip channel: lists of tasks and progress updates were regularly posted, and volunteers could claim the tasks and report on their progress.
The correctness of the code was ensured by continuous integration scripts on github, which automatically check that the code compiles without errors, and by manual review of the proposed additions.
The blueprint continued to evolve during the formalization, as new lemmas were added (once because a small gap was identified in the blueprint proof), definitions were changed to better fit the needs of writing in \Lean, and small mistakes were corrected (particularly in the exact value of constants used in the Kolmogorov-Chentsov proof).

\subsection{Organization of the paper}

Section~\ref{S:overview} gives an overview of the different theorems and definitions we formalized in order to construct Brownian motion.
Section~\ref{S:RV-formalism} discusses the random variable formalism we developed in \mathlib as part of this project.
Section~\ref{S:formalization} describes in more detail our formalization of the main theorems: Carathéodory and Kolmogorov extension theorems, definitions and properties of Gaussian measures, and Kolmogorov-Chentsov continuity theorem.
Finally, Section~\ref{S:BM} contains the final construction of Brownian motion and the Wiener measure on $\R_+$, and a discussion of the properties of Brownian motion that we formalized.

\section{Formalization overview}\label{S:overview}
This section gives an overview of the different theorems and definitions we formalized in order to construct Brownian motion. Subsequent sections will go into more detail for each of these steps.
We refer to the implementation using links, like \mllink{Probability/Distributions/Gaussian/Basic}{L43-L46}{ProbabilityTheory.IsGaussian}, which point either to \mathlib or to \href{https://github.com/RemyDegenne/brownian-motion}{our repository}.
Our goal is to eventually integrate all of our code into \mathlib, and that process is ongoing.
All links point to commit {\tt 82d756ea359cfa0f2c10bc00b8a997822c3b32ea} for \mathlib, and to commit {\tt 6a7a3141c6bc217484f635f3ab861316a37c0362} for our repository, which correspond its state on December 1, 2025.

Let us start with notation. Any stochastic process has an index set $T$ (usually called {\em time}, usually uncountably infinite), and a state space.
While the general use of the word {\em time} implies a (total) order, the index set $T$ does not need to be ordered for the general construction of stochastic processes. Only for the definition of the covariance matrix in Section~\ref{ss:char}, we need the order relation. In addition, $T$ is assumed to be metric in Section~\ref{ss:kolchen}. For the state space, since we use dependent type theory, the state space may depend on the time, i.e.\ the state space at time $t$ is $E_t$.
This is interesting, since virtually no standard textbook treats the case of a time-dependent state space. We only found a few examples of time-dependent stochastic processes in the literature (one of which is \cite{arnaudon2008brownian}).
Extending the usual hypothesis on the state space, we require that $E_t$ is a complete and second-countable metric space, for all $t \in T$. Probability measures are denoted $P_.$.

\subsection{The Carathéodory and Kolmogorov Extension Theorems}
The usual approach to construct (the distribution of) a stochastic process works as follows: describe properties of the distribution of the stochastic process $P_J$ at some arbitrary but finite number of times $J = \{t_1,...,t_n\} \subseteq T$.
The resulting family of probability measures $(P_J)_{J \subseteq T \text{ finite}}$ has to be {\em projective} in the sense that the projection of $P_J$ to $H\subseteq J$ has to be equal to $P_H$. In other words, when describing the distribution of the stochastic process at all times in $J$, and then forgetting all properties for times in $J\setminus H$, results in the description of properties at times in $H$. One may then ask if this already gives a complete description of the process for all times.

It is the achievement of Kolmogorov that these finite-dimensional distributions in fact provide a unique description of the distribution of a stochastic process, even if $T$ is uncountable, as long as the underlying family of state spaces $(E_t)_{t\in T}$ is nice enough (Polish, i.e.\ a second-countable topological space which can be metrized by a complete metric, for example) \cite{kolmogoroff1933grundbegriffe}.
The resulting measure is defined on the product-$\sigma$-field $\mathcal F :=\bigotimes_{t\in T} \mathcal B(E_t)$ (where $\mathcal B(E_t)$ is the Borel $\sigma$-algebra on $E_t$). Here, $\mathcal F$ is generated by finite projections and hence any element of $\mathcal F$ may only depend on at most countably many $t\in T$, making this a rather coarse $\sigma$-algebra. (In particular, note that this is not the Borel $\sigma$-algebra of the product topology for uncountable $T$.)

The Kolmogorov extension theorem is on the interface between measure theory and probability theory. Here, we rely on a decent amount of formalized mathematics in the measure-theory part of \mathlib (outer measures, above all), while not requiring any specific previous formalization of probability theory. (In fact, most of our results are formulated in terms of finite rather than probability measures.)

We are going to formulate the main result in a modern fashion, as e.g.\ found in Theorem 2.2 of \cite{rao1971projective}, Theorem 7.7.1 of Volume~2 of \cite{bogachev2007measure}, Theorem 15.26 of \cite{guide2006infinite}, or \cite{border1998expository}. Note that these formulations split general assumptions on the underlying space(s) (e.g.\ a metric property) from the property which is needed in the proof (inner regularity with respect to compact sets). Other -- highly readable -- references such as \cite{Billingsley1995} state the extension theorem only in special cases such as $E_t = \mathbb R$ for all $t$.

\subsection{Gaussian Measures and characteristic functions}
\label{ss:char}
Our goal in this subsection is to define the finite-dimensional distributions of Brownian motion. For times $t := (t_1, ..., t_n)$, this is given by the multi-dimensional Gaussian distribution $N(0, C_t)$, where $0$ is the vector of expectations, and $C_t$ is the covariance matrix, given by $(C_t)_{ij} = t_i \wedge t_j$. In order to do so, we rely on two tools that were already in \mathlib, namely the implementation of the one-dimensional normal distribution (\mllink{Probability/Distributions/Gaussian/Real}{L198-L201}{ProbabilityTheory.gaussianReal}) and characteristic functions of probability (or finite) measures.

For any probability measure $P$ on $\mathbb R^n$, its characteristic function is given by $\psi: t \mapsto \int e^{it^\top x} P(dx)$, where the integral takes values in $\mathbb C$. We use the fact that $\psi$ characterizes $P$ uniquely (\mllink{MeasureTheory/Measure/CharacteristicFunction}{L235-L244}{MeasureTheory.Measure.ext_of_charFun}).
For the standard normal distribution $N(0,1)$ this is $\psi_{N(0,1)}(t) = \exp(-t^2/2)$.
It so happens that a probability measure $P$ on $\mathbb R^n$ is Gaussian if and only if its characteristic function is given by $\psi(t) = \exp(im - t^\top C t / 2)$, where $m$ is the vector of expectations and $C$ is the covariance matrix.
We can define the $n$-fold product measure $N(0, I_n) := N(0, 1)^{\otimes n}$, where $I_n$ is the unit matrix, with characteristic function $\psi_{N(0,I_n)}(t) = \exp\big(-\tfrac 12 t^\top I_n t\big)$.
It is easy to see that the covariance matrix of a measure must be positive semidefinite. Conversely, for any such $C$, there is $A$ with $C = A^\top A$ and we can take the image of $N(0,I_n)$ under the map $f : x\mapsto Ax$. It has the characteristic function
\begin{align} \label{eq:gausslin}
\psi_{f_\ast N(0,I_n)}(t) = \int e^{it^\top x} f_\ast N(0,I_n) dx = \int e^{it^\top A x} N(0,I_n) dx = \exp\big( - \tfrac 12 t^\top C t\big).
\end{align}
So, we can define the finite dimensional distributions of Brownian motion given that we show that $C \in \R^n$ with entries $C_{ij} = t_i \wedge t_j$ is positive semidefinite.
We can also show with characteristic functions that the family of finite-dimensional distributions of Brownian motion is projective, as needed for the Kolmogorov extension theorem.

\subsection{The Kolmogorov-Chentsov Theorem}
\label{ss:kolchen}

Brownian motion is a process with continuous paths. To show the existence of a modification of a stochastic process with continuous paths, we use the Kolmogorov-Chentsov Theorem.
Let us describe briefly that theorem in a simple case here, before we dive deeper in the next sections.
Let $X$ be a stochastic process on the index set $T = [0,1]$ with values in a metric space $E$ with distance $d_E$.
The result states that if we assume we find $p>0, q>1, M > 0$ satisfying
\begin{align}
\label{eq:cs}
  \mathbf E[d_E(X_s, X_t)^p] \leq M|t-s|^q, \qquad 0\leq s,t\leq 1\:,
\end{align}
then there exists $Y = (Y_t)_{t\in [0,1]}$ with $\mathbf P(X_t = Y_t) = 1$ for all $t\in [0,1]$ and $Y$ has almost surely Hölder continuous paths with coefficient $\gamma$ for any $\gamma < \tfrac {q-1} p$.

The most commonly given proof of this statement is a follows.
Set $D_n := \{k/2^n: k=0,...,2^n\}$ and $D := \bigcup_{n\in\mathbb N} D_n$, the set of dyadics. Start by showing summability of $\mathbf P\Big( \sup_{s,t\in D_n, |t-s| = 2^{-n}} d_E(X_s, X_t) \geq 2^{-\gamma n} \Big)$ using \eqref{eq:cs}.
By the Borel-Cantelli Lemma, this shows that $\sup_{s,t\in D_n, |t-s| = 2^{-n}} d_E(X_s, X_t) \leq 2^{-\gamma n}$ holds for $n$ large enough.
From this, we see that $X$ is locally Hölder-$\gamma$ continuous on $D$.
From here, define some Hölder continuous $Y$ coinciding with $X$ on $D$.
Finally, fix $t \in [0,1]$ and $t_1, t_2,...\in D$ with $t_n \xrightarrow{n\to\infty} t$.
Using \eqref{eq:cs}, we see that $X_{t_n} \xrightarrow{n\to\infty} X_t$ in probability, as well as $Y_{t_n} \xrightarrow{n\to\infty} Y_t$ almost surely due to continuity of $Y$.
Therefore, $X_t = Y_t$ almost surely.

We will use a more general version of this statement, replacing $T = [0,1]$ by a metric space with a property restricting the number of balls needed to cover $T$. The version we formalize is based on
the recent work of~\cite{kratschmer2023kolmogorov}.

\subsection{Construction of Brownian motion and Wiener measure on $\R_+$}

In order to finally construct Brownian motion, we need to put everything together.
We have constructed a projective family of finite-dimensional distributions given by $N(0,C_t)$ for $t = (t_1,...,t_n)$ with $(C_t)_{ij} = t_i \wedge t_j$.
We can then invoke the Kolmogorov extension theorem to obtain a probability measure on $\mathbb{R}^{\mathbb{R}_+}$ and define a stochastic process $X = (X_t)_{t\geq 0}$ with that distribution.
Finally, the assumption in the Kolmogorov-Chentsov Theorem can be verified by computing moments of Gaussian distributions.
For example,
\[ \mathbf E[|X_t - X_s|^4] = \mathbf E[X_{t-s}^4] = (t-s)^2 \mathbf E[X_{1}^4] < \infty.\]
So, we obtain from that theorem a process with continuous paths and the correct finite dimensional distributions, which we call Brownian motion.
This also gives a distribution on the continuous functions $\mathcal C(\mathbb R_+, \mathbb R)$, which is the Wiener measure.

\section{Random variable formalism}
\label{S:RV-formalism}

Before describing the formalization of the theorems leading to the Brownian motion construction, we discuss how this project led to an interesting development of the random variable formalism in \mathlib. Most content about probability theory in the library focuses on measures, as opposed to random variables.
The library contains several notions related to random variables, such as independence, conditional distributions or sub-Gaussianity. However it remains quite cumbersome to deal with random variables only without ever going back to measure theory. This strongly contrasts with probability theory in paper mathematics, where we set a probability space at the beginning and then do all the reasoning on random variables.

If we had stuck to measure theory, we would just have built a measure on the space of continuous functions and proved properties about the canonical process. But in paper maths, we more generally consider a Brownian motion over a general probability space, discuss its finite dimensional laws, its independence with respect to other processes, or its invariance properties. All these concepts can only be discussed in a random variable framework. This led us to introduce many definitions and results to be able to manipulate these random variables in a way that is similar to what is done in paper maths.

\paragraph{Law of a random variable}

The law of a random variable is a key concept in probability theory. Given a probability space $(\Omega, \mathcal{F}, P)$ and a random variable $X : \Omega \to E$, the law of $X$ is the pushforward $X_*P$ of the measure $P$ by the map $X$. In \Lean this is spelt as \lean{P.map X}. This however already comes with an issue. When one writes that $X$ is a random variable, it is implied that the map is required to be measurable. In \Lean however this has to be added as a separate hypothesis.
This can quickly get annoying.
For instance, say $X$ has law $\mu$ under $P$ and one wants to compute $P[f(X)]$ for some measurable function $f$ to show that it is equal to $\int f(x) \mu(\mathrm dx)$. Then one would have to rewrite saying that $P[f(X)]$ is actually $\int (f \circ X)(\omega) P(\mathrm d\omega) = \int f(x) (X_*P)(\mathrm dx) = \int f(x) \mu(\mathrm dx)$. The first equality in particular is not as straightforward as one would like formally, first because \Lean has to recognize the integrand as being of the form $(f \circ X)(\omega)$, which does not always work when $f$ or $X$ has a complicated form,
and then because one must also provide a proof that $X$ is measurable.
When this occurs many times in a proof, it becomes increasingly tedious.

To tackle this issue, we introduced a new predicate \lean{HasLaw X µ P}. This is a bundled statement that contains both the fact that \lean{P.map X = µ} and the fact that $X$ is measurable. The computation above can the be contained in only one lemma, namely \mllink{Probability/HasLaw}{L100-L104}{ProbabilityTheory.HasLaw.integral_comp}. But it goes further than that, because capturing this concept into a predicate also allows one to manipulate random variables without ever referring to the pushforward of a measure, but only manipulating the laws of the random variables. For instance, if the function $f$ maps the measure $\mu$ to a measure $\nu$, we have the lemma \mllink{Probability/HasLaw}{L72-L77}{ProbabilityTheory.HasLaw.comp} which states \lean{HasLaw (f ∘ X) ν P}. We can likewise state that the law of the sum of two independent random variables is the convolution of their laws. This gives a nice framework for computations with random variables.

\paragraph{Gaussian random variables}

Part of the work in this project is to prove results about Gaussian measures. However in the end what we want is to manipulate Gaussian \emph{random variables}, and it would be inconvenient to have to write \lean{IsGaussian (P.map X)} and deal with pushforward measures once again. This is why we introduced another predicate on random variables: \lean{HasGaussianLaw X P}.
This is precisely the same statement, but written in a more natural way which avoids mentioning pushforward of measures. We can then prove for instance that if \lean{HasGaussianLaw X P} and and we consider \lean{L} a continuous linear map, then \lean{HasGaussianLaw (L ∘ X) P} (\bmlink{Auxiliary/HasLaw}{L118-L122}{ProbabilityTheory.HasGaussianLaw.map}).
Specializing this statement to different values for \lean{L} allows for instance to prove that if $(X_i)_{1 \le i \le n}$ is a Gaussian vector, i.e.\ \lean{HasGaussianLaw (fun ω ↦ (X · ω)) P}, then each $X_i$ is also Gaussian, i.e.\ \lean{HasGaussianLaw (X i) P}.

One of the great properties of Gaussian random variables is their relation to independence. Proving that two random variables are independent can be quite involved sometimes. But when it comes to Gaussian random variables, we have the following result. If $X$ and $Y$ are two real random variables such that the pair $(X, Y)$ is Gaussian and such that $\mathbf{cov}(X, Y) = 0$, then $X$ and $Y$ are independent. In other words, independence can be checked by only computing a covariance. This is not only true for real variables and can be generalized. The most general form we formalized is stated below.

\begin{lemma}[\bmlink{Auxiliary/HasGaussianLaw}{L97-L103}{HasGaussianLaw.iIndepFun_of_cov}]\label{lem:gauss-indep}
  Consider $(E_i)_{1 \le i \le n}$ a family of second countable Banach spaces, and a family of random variables $X_i : \Omega \to E_i$. Assume that the joint random variable $(X_i)_{1 \le i \le n}$ is Gaussian. Then to check that the $X_i$ are mutually independent, it is enough to check that for any $1 \le i \ne j \le n$ and any $L_1 \in E_i'$, $L_2 \in E_j'$, we have $\mathbf{cov}(L_1(X_i), L_2(X_j)) = 0$.
\end{lemma}

This statement follows actually quite easily from the general results we formalized about characteristic functions of Gaussian measures (see Section~\ref{S:gaussian}).

\paragraph{Gaussian processes}

We end this section by presenting our formalization of Gaussian processes. These are of course crucial if we want to talk about the Brownian motion beyond the measure theoretic aspects, as the Brownian motion is itself a Gaussian process.

\begin{definition}[\bmlink{Gaussian/GaussianProcess}{L27-L30}{ProbabilityTheory.IsGaussianProcess}]\label{def:guass-proc}
  Consider $E$ a Banach space and $(X_t)_{t \in T}$ a stochastic process taking values in $E$. Then $X$ is a \emph{Gaussian process} if for all $t_1, ..., t_n \in T$, the variable $(X_{t_1}, ..., X_{t_n})$ is Gaussian.
\end{definition}

This can be stated in a very straightforward way thanks to the \lean{HasGaussianLaw} predicate, and is bundled into the predicate \lean{IsGaussianProcess X P}. Once again this illustrates our goal of relying on statements about random variables. Gaussian processes enjoy just as nice properties as Gaussian random variables, in particular with respect to independence. We have for example the following lemma.

\begin{lemma}[\bmlink{Gaussian/GaussianProcess}{L201-L204}{IsGaussianProcess.indepFun''}]\label{lem:indepGaussProc}
  Let $(X_s)_{s \in S}$ and $(Y_t)_{t \in T}$ be two stochastic processes taking values in $E$. Assume that the joint process $((X_s)_{s \in S}, (Y_t)_{t \in T})$ is Gaussian. Then to show that $X$ and $Y$ are independent, it is enough to prove that for any $s \in S$ and $t \in T$, $\mathbf{cov}(X_s, Y_t) = 0$.
\end{lemma}

This lemma reduces the independence of stochastic processes to a very simple computation, and will prove to be a blessing to show that a Brownian motion enjoys the weak Markov property. It is the result of a combination of Lemma~\ref{lem:gauss-indep} and of some results about independence of processes which were added to \mathlib for this purpose (notably \mllink{Probability/Independence/Process}{L164-172}{ProbabilityTheory.IndepFun.process_indepFun_process}).

The statement requires $((X_s)_{s \in S}, (Y_t)_{t \in T})$ to be Gaussian. However, written like that, it is not technically a process, but a pair of processes. To state it cleanly in \Lean we consider the disjoint sum $S \oplus T$. This is the set where each element is either an element of $S$, or an element of $T$, both options being exclusive. It corresponds to \texttt{Sum} in \Lean. We can then consider the process
\begin{align*}
  X \oplus Y : S \oplus T & \to E \\
  S \ni s & \mapsto X_s. \\
  T \ni t & \mapsto Y_t
\end{align*}

We formalized more results about Gaussian processes that we present in the next section to see how they were applied to prove properties of the Brownian motion. This idea of defining predicates on random variables rather than measures proved to be really nice to work with, landing a formalization process of basic properties of Brownian motion which felt quite close to the paper proofs.
This demonstrates that it is possible to mimic basic reasoning on random variables in \Lean without suffering of some obstacles due to formalization.

\section{The formalization}\label{S:formalization}

We describe in more detail our formalization of the main results outlined in Section~3.
After a few preliminaries about stochastic processes in Section~\ref{S:stochastic_processes} we first present the Carathéodory and Kolmogorov Extension Theorems in Section~\ref{S:extension}, then Gaussian measures and characteristic functions in Section~\ref{S:gaussian}, and finally the the Kolmogorov-Chentsov Theorem in Section~\ref{S:continuity}.

\subsection{Stochastic processes}\label{S:stochastic_processes}

Let $\Omega$ be a measurable space, on which we define a probability measure $P$ (with expectation denoted by $\mathbf{E}$).
Let $T$ be a set, and for each $t \in T$, let $E_t$ be a measurable space.
A stochastic process indexed by $T$ with values in $(E_t)_{t \in T}$ is a collection $(X_t)_{t \in T}$ of random variables $X_t : \Omega \to E_t$, which satisfies a measurability condition (see further down).
We will often take $E_t = E$ for all $t \in T$ for some fixed measurable space $E$, but we first prove general results about processes with values in different spaces.
There is no \Lean definition for a stochastic process: we use a function $X : (t : T) \to \Omega \to E_t$ together with appropriate measurability conditions.

The goal of our formalization of Brownian motion is to define a process with a given law and continuous paths.
The law of a stochastic process is the pushforward of the measure $P$ by the map $\Omega \to \prod_{t \in T} E_t$ defined by $\omega \mapsto (X_t(\omega))_{t \in T}$.
The paths of the process are the maps $(t : T) \to E_t$ defined by $t \mapsto X_t(\omega)$ for each $\omega \in \Omega$.

A process $Y$ is said to be a modification of $X$ if for all $t \in T$, $X_t = Y_t$ almost surely.
$X$ and $Y$ are said to be indistinguishable if almost surely, $X_t = Y_t$ for all $t \in T$ simultaneously.
Being indistinguishable implies being a modification, which in turn implies equality of laws (that last part is proved through finite dimensional distributions, see below).
We will use repeatedly that if $T$ is a second-countable topological space, $E$ is Hausdorff and the paths of $X$ and $Y$ are continuous almost surely, then being a modification implies being indistinguishable (\bmlink{Gaussian/StochasticProcesses}{L14-L18}{indistinguishable_of_modification}).

\paragraph{Measurability}

Working with stochastic processes requires a form of measurability.
First, we say that $X$ is measurable if $X_t$ is measurable for all $t \in T$.
This is equivalent to requiring that the map $\Omega \to \prod_{t \in T} E_t$ defined by $\omega \mapsto (X_t(\omega))_{t \in T}$ is measurable for the product $\sigma$-algebra on $\prod_{t \in T} E_t$.

We say that $X$ is almost everywhere measurable if there exists a measurable process $Y$ such that almost surely, $X_t = Y_t$ for all $t \in T$.
That is, there exists a measurable process $Y$ which is indistinguishable from $X$.
Contrary to the case of measurability, this is strictly stronger than requiring that for all $t \in T$, $X_t$ is almost everywhere equal to a measurable function (unless for example $T$ is countable or $X$ has continuous paths, in which case the two notions coincide).
In the following, a \emph{stochastic process} will always mean an almost everywhere measurable collection $(X_t)_{t \in T}$.

\mathlib uses extensively almost everywhere measurable random variables, more so than measurable ones, because almost everywhere measurability is a much more flexible notion, which does not depend on the behavior of the function on sets to which the measure gives zero mass.
\cite{gouezel2022formalization} contains a detailed discussion of the advantages of almost everywhere measurability when it comes to defining integrals in the context of calculus.

\paragraph{Finite dimensional distributions}

\sloppy The finite dimensional distributions of a stochastic process $(X_t)_{t \in T}$ are the pushforwards of $P$ by the maps $\Omega \to \prod_{i=1}^n E_{t_i}$ defined by $\omega \mapsto (X_{t_1}(\omega),...,X_{t_n}(\omega))$ for all $n \in \mathbb{N}$ and all choices of $t_1,...,t_n \in T$.
Their importance comes from the fact that two stochastic processes have the same law if and only if they have the same finite-dimensional distributions, which we prove in \mllink{Probability/Process/FiniteDimensionalLaws}{L75-L81}{ProbabilityTheory.identDistrib_iff_forall_finset_identDistrib} by using the uniqueness of the projective limit of a family of measures described further down.
Indeed the law of a process is the projective limit of its finite-dimensional distributions.
Our goal in this project is thus to build a process with given finite-dimensional distributions and continuous paths.

\subsection{The Carathéodory and Kolmogorov Extension Theorems}
\label{S:extension}
We describe our implementation of the Kolmogorov Extension Theorem, which relies on basic notions from topology such as metric spaces and the Borel $\sigma$-algebra.
We begin by defining set systems we will need, which correspond to the definitions \mllink{MeasureTheory/SetSemiring}{L61-L67}{MeasureTheory.IsSetSemiring} and \mllink{MeasureTheory/SetSemiring}{L453-L458}{MeasureTheory.IsSetRing} (which we added to \mathlib).

\begin{definition}[Semi-ring, ring]\label{def:semi}
  Let $E$ be some set. We call a set of sets $\mathcal H \subseteq 2^E$ a
  \emph{semi-ring}, if it is (i) a $\pi$-system (i.e.\ closed under
  $\cap$) and (ii) for all $A, B \in\mathcal H$ there is\footnote{We write $A \subseteq_f B$ if $A$ is a finite subset of $B$.} $\mathcal K
  \subseteq_f \mathcal H$ with\footnote{We write $A\uplus B$ for
    $A\cup B$ if $A\cap B=\emptyset$.}  $B\setminus A = \biguplus_{K
    \in \mathcal K} K$.  \\ We call $\mathcal H \subseteq 2^E$ a
  \emph{ring}, if it is closed under $\cup$ and under set-differences.
\end{definition}

\noindent
Any ring is a semi-ring since $A\cap B = A \setminus (A \setminus B)$, i.e.\ every ring is a $\pi$-system, and the second condition is satisfied with $\mathcal{K}=\{B \setminus A\}$.
Let us state two important lemmas on semi-rings.

\begin{lemma}[\mllink{MeasureTheory/SetSemiring}{L221-L228}{MeasureTheory.IsSetSemiring.disjointOfDiffUnion}]\label{l1}
  Let $\mathcal H$ be a semi-ring, $\mathcal I \subseteq_f \mathcal H$,
  $A \in \mathcal H$. Then, there is $\mathcal K \subseteq_f \mathcal
  H$ such that $\mathcal K$ contains pairwise disjoint sets and $A
  \setminus \bigcup_{I \in \mathcal I} I = \biguplus_{K\in \mathcal K}
  K$.
\end{lemma}

\begin{lemma}[\mllink{MeasureTheory/SetSemiring}{L405-L411}{MeasureTheory.IsSetSemiring.disjointOfUnion}]\label{l2}
  Let $\mathcal H$ be a semi-ring and $A_1,...,A_m \in \mathcal
  H$. Then, there are $\mathcal K_1,...,\mathcal K_m \subseteq_f
  \mathcal H$ disjoint such that $\bigcup_{n=1}^m \mathcal K_n$
  contains disjoint sets and $\bigcup_{m=1}^n A_m = \biguplus_{m=1}^n
  \biguplus_{K \in \mathcal K_n} K$.
\end{lemma}

Given an additive content $m : \mathcal H \to [0,\infty]$ for some semi-ring $\mathcal H$ (see Definition~\ref{def:content} as well as \mllink{MeasureTheory/Measure/AddContent}{L66-L74}{MeasureTheory.AddContent}), the goal of Carathéodory's extension theorem is to define a measure $\mu : \sigma(\mathcal H) \to [0,\infty]$ extending $m$ to the $\sigma$-algebra $\sigma(\mathcal H)$ generated by $\mathcal H$. More precisely, the Carathéodory extension gives a measure on an even larger $\sigma$-algebra; see Theorem~\ref{T:cara} in conjunction with Theorem~\ref{T:masseind} below. We will follow this abstract construction, and start by stating some basic concepts; see \mllink{MeasureTheory/Measure/MeasureSpaceDef}{L71-L81}{MeasureTheory.Measure}, \mllink{MeasureTheory/OuterMeasure/Defs}{L50-L57}{MeasureTheory.OuterMeasure}.

\begin{definition}\label{def:content}
  For some set $E$, let $\mathcal H \subseteq 2^E$
  and call any $m : \mathcal H \to [0,\infty]$ a content.
  \begin{enumerate}
  \item $m$ is called additive if for $\mathcal K \subseteq_f \mathcal H$ pairwise disjoint and $\biguplus_{K \in \mathcal K} K \in \mathcal H$, we have $m \Big(\biguplus_{K \in \mathcal K} K \Big) = \sum_{K \in \mathcal K} m(K)$. If the same holds for\footnote {We write $A \subseteq_c B$ if $A$ is a countable subset of $B$.}$\mathcal K \subseteq_c \mathcal H$ pairwise disjoint, we say that $m$ is $\sigma$-additive.
  \item The set-function $m$ is called sub-additive if for $\mathcal K \subseteq_f \mathcal H$ and $\bigcup_{K \in \mathcal K} K \in \mathcal H$, we have $m \Big(\bigcup_{K \in \mathcal K} K \Big) \leq \sum_{K \in \mathcal K} m(K)$. (Note that elements of $\mathcal K$ need not be disjoint.) Here, $\sigma$-sub-additivity is defined in the obvious way using $\mathcal K\subseteq_c \mathcal H$.
  \item If $m(A) \leq m(B)$ for $A\subseteq B$ and $A,B\in\mathcal H$, we say that $m$ is monotone.
  \item If $\mathcal H$ is a $\sigma$-algebra and $m$ is $\sigma$-additive with $m(\emptyset) = 0$, we call
    $m$ a measure.
  \item If $\mathcal H = 2^E$, $m$ is monotone and $\sigma$-sub-additive with $m(\emptyset)=0$, we call $m$ an outer measure.
  \end{enumerate}
\end{definition}

\subsubsection{Carathéodory's Extension Theorem}
An additive, $\sigma$-sub-additive content on a semi-ring induces an outer measure by approximating sets from above (see \mllink{MeasureTheory/OuterMeasure/Induced}{L293-L295}{MeasureTheory.inducedOuterMeasure_eq}):

\begin{proposition}[Outer measure induced by a set function on a semi-ring] 
  Let \label{P:auss} $\mathcal H$ be a semi-ring and $m: \mathcal H\to\mathbb R_+$ additive and $\sigma$-sub-additive. For $A\subseteq E$ let
  $$ \mu(A) := \inf_{\mathcal G \in \mathcal U(A)} \sum_{G\in\mathcal G} m(G)$$ where
  $$ \mathcal U(A) := \big\{\mathcal G \subseteq_c \mathcal H, A\subseteq \bigcup_{G\in\mathcal G} G\big\}$$
  is the set of countable coverings of $A$. Then, $\mu$ is an outer measure.
\end{proposition}

Moreover, we call a set Carathéodory with respect to some outer measure, if it is measurable in the following sense (see \mllink{MeasureTheory/OuterMeasure/Caratheodory}{L53-L56}{MeasureTheory.OuterMeasure.IsCaratheodory}):

\begin{theorem}[\boldmath $\mu$-measurable sets are a
    $\sigma$-algebra]\label{T:cara} Let $\mu$ be an outer measure on
  $E$ and $\mathcal F$ the set of $\mu$-measurable sets,
  i.e.\ the set of sets $A$ satisfying
  \begin{align*}
    \mu(B) = \mu(B\cap A) + \mu(B\cap A^c), \qquad B \subseteq E.
  \end{align*}
  Then, $\mathcal F$ is a $\sigma$-Algebra and $\mu|_{\mathcal F}$ is
  a measure.
\end{theorem}

In the construction above, since all sets in the semi-ring $\mathcal H$ are Carathéodory (see \mllink{MeasureTheory/OuterMeasure/OfAddContent}{L86-L88}{MeasureTheory.AddContent.isCaratheodory_ofFunction_of_mem}), the $\sigma$-algebra from Theorem~\ref{T:cara} is at least as large as the $\sigma$-algebra generated by $\mathcal H$. The resulting measure is \mllink{MeasureTheory/OuterMeasure/OfAddContent}{L134-L138}{MeasureTheory.AddContent.measureCaratheodory} and agrees with the additive content on $\mathcal H$ (see \mllink{MeasureTheory/OuterMeasure/OfAddContent}{L171-L176}{MeasureTheory.AddContent.measure_eq}). Let us summarize this:

\begin{theorem}[Carathéodory extension]\label{T:masseind}
  Let $\mathcal H$ be a semi-ring and $m: \mathcal H\to\mathbb R_+$
  $\sigma$-finite and $\sigma$-additive. Furthermore, let $\mu$ be the
  induced outer measure from Proposition~\ref{P:auss} and $\mathcal F$
  the $\sigma$-algebra from Theorem~\ref{T:cara}. Then,
  $\sigma(\mathcal H)\subseteq\mathcal F$ and $\mu$ coincides with $m$
  on $\mathcal H$.
\end{theorem}

In order to apply Theorem~\ref{T:masseind}, we need to show that an additive content is $\sigma$-additive.
We do this on the level of rings rather than semi-rings. Note that we can extend an additive content on a semi-ring to an additive content on the ring generated by the semi-ring.

\begin{lemma}
Let $\mathcal H$ be a semi-ring. Then,
\begin{align*}
\mathcal R := \Big\{ \biguplus_{j \in J} A_j : J \text{ finite, } A_j \in \mathcal H, (A_j)_{j \in J} \text{ disjoint}\Big\} \supseteq \mathcal H
\end{align*}
is a ring. Moreover, if $m$ is an additive content on $\mathcal H$, $m' : \mathcal R \to \mathbb R_+$, defined, for disjoint $(A_j)_{j\in J}$ in $\mathcal H$
$$ m'\Big( \biguplus_{j \in J} A_j \Big) := \sum_{j\in J} m(A_j),$$
extends $m$ to $\mathcal R$.
\end{lemma}

Actually, $\sigma$-additivity and $\sigma$-sub-additivity are close relatives.
We first have the result \mllink{MeasureTheory/Measure/AddContent}{L388-L397}{MeasureTheory.addContent_iUnion_eq_sum_of_tendsto_zero}, which states that $\sigma$-additivity is implied by continuity of the content at the empty set, and \mllink{MeasureTheory/Measure/AddContent}{L437-L441}{MeasureTheory.isSigmaSubadditive_of_addContent_iUnion_eq_tsum}, which states that -- if the additive content is defined on a ring -- it is $\sigma$-sub-additive if it is $\sigma$-additive.
Therefore, a convenient way to show $\sigma$-additivity of a content (in order for example to apply Carathéodory's extension theorem) is to show continuity in $\emptyset$.
For this, we require the notion of inner regularity with respect to a compact system.

We recall the definition of inner regularity of set functions, which was already in \mathlib.

\begin{definition}[Inner regularity \mllink{MeasureTheory/Measure/Regular}{L203-L210}{MeasureTheory.Measure.InnerRegularWRT}] \label{def:innerreg}
  Let $E$ be some set, equipped with a topology, and $m$ be a
  set-function on some $\mathcal H \subseteq 2^E$.
  \begin{enumerate}
  \item Let $p, q : 2^E \to \{\text{true, false}\}$. Then, $m$ is
    called inner regular with respect to $p$ and $q$, if
    $$ m(A) = \sup\{m(F) : p(F) = \text{true}, F \subseteq A\}$$ for
    all $A \in \mathcal H$ with $q(A) = \text{true}$.
  \item If $q(A) = \text{true}$ iff $A$ is measurable, we write simply ``inner regular with respect to $p$'' (omitting ``and $q$''). If $p(A) = \text{true}$ iff $A$ is closed (compact,
    closed and compact), we say that $m$ is inner regular with respect
    to closed (compact, compact and closed) sets.
  \end{enumerate}
\end{definition}

For the next result, recall that for compact sets $C_1, C_2,...$ with
$\bigcap_{n=1}^\infty C_n = \emptyset$, there is some $N$ with
$\bigcap_{n=1}^N C_n = \emptyset$. More generally, compact sets form a
compact system, which is defined as follows.

\begin{definition}[Compact system, \kelink{CompactSystem}{L17-L18}{IsCompactSystem}]
  Let $\mathcal C \subseteq 2^E$. If, for all $C_1, C_2,...$ with
  $\bigcap_{n=1}^\infty C_n = \emptyset$, there is some $N$ with
  $\bigcap_{n=1}^N C_n = \emptyset$, we call $\mathcal C$ a {\em
    compact system}.
\end{definition}

\noindent
Such compact systems are important since they allow for a proof of continuity at $\emptyset$ (hence $\sigma$-additivity) of a content on a ring, which is the missing piece for applying the Carathéodory Theorem in the proof of the Kolmogorov extension theorem; see \kelink{RegularContent}{L17-L21}{MeasureTheory.tendsto_zero_of_regular_addContent}.

\begin{lemma}\label{l:stetigcompact}
  Let $\mu$ be an additive set function on a ring $\mathcal R$, which contains the compact system $\mathcal C$. If $\mu$ is inner regular with respect to $\mathcal C$, then $\mu$ is continuous at $\emptyset$.
\end{lemma}

\subsubsection{Kolmogorov's Extension Theorem}
We aim to apply Theorems~\ref{T:cara} and~\ref{T:masseind} on product spaces. The resulting Kolmogorov's extension theorem is a statement about extending a set-function on a product space to a (finite) measure, where the product space can come with an arbitrary index set.  The next definition covers the important concept of a projective family of measures. In short, we define measures on any finite subset of indices in a consistent way; see \mllink{MeasureTheory/Constructions/Projective}{L43-L47}{MeasureTheory.IsProjectiveMeasureFamily}.

\begin{definition}[Projective family and projective limit] \mbox{}
  \begin{enumerate}
  \item Let $T$ be some (index) set and $(E_i)_{i\in T}$ a
    family of sets, such that $\mathcal F_i$ is a $\sigma$-algebra on $E_i, i \in T$. For $J\subseteq T$, we denote $E_J :=
    \prod_{j \in J} E_j$ and $\pi_J : E_T \to E_J$
    the projection. For $H\subseteq J \subseteq T$, we write
    $\pi_H^J$ for the projection $E_J \to E_H$.
  \item Let $\mathcal F_i$ be a $\sigma$-algebra on $E_i$,
    $i\in T$. For $J\subseteq_f T$, let $\mathcal F_J$ be the
    (finite) product-$\sigma$-algebra on $E_J$, and $\mathcal F_T$ be
    the $\sigma$-algebra generated by cylinder sets
    $\{\pi_J^{-1} A: J \subseteq_f T, A \in
    \mathcal F_J\}$.
  \item A family $(P_J)_{J\subseteq_f I}$, where $P_J$ is a finite
    measure on $\mathcal F_J$, is called projective if
    $$P_H = (\pi_H^{J})_\ast P_J$$ for all
    $H\subseteq J \subseteq_f I$. (Recall that $A \mapsto
    (\pi_H^{J})_\ast P_J(A) := P_J((\pi_H^{J})^{-1}A)$
    is called the pushforward measure of $P_J$ under $\pi_H^{J}$.)
  \item If, for some projective family $(P_J)_{J\subseteq_f T}$,
    there is a finite measure $P_T$ on $\mathcal F_T$ with
    $P_J = (\pi_J)_\ast P_T$ for all $J\subseteq_f T$, then we
    call $P_T$ a projective limit of $(P_J)_{J\subseteq_f T}$.
  \end{enumerate}
\end{definition}

In order to apply Theorem~\ref{T:masseind} and show the extension theorem, we only need to show that $\{\pi_J^{-1} C: C \in E_J \text{ compact and closed}\}$ is a compact system. Note that compact sets are closed in Hausdorff spaces, but we do not have this property since we are working with pseudo-metric spaces, which do not have this property.\footnote{Recall, that a pseudo-metric $r(.,.)$ can have $r(x,y) = 0$ for $x\neq y$, and therefore is not a $T2$-space.}
Since Lemma~\ref{l:stetigcompact} gives the $\sigma$-additivity of an additive content \mllink{MeasureTheory/Constructions/ProjectiveFamilyContent}{L114-L120}{MeasureTheory.projectiveFamilyContent}, which is defined through the projective family $P$, we have a proof of the Kolmogorov Extension Theorem \kelink{KolmogorovExtension}{L95-L99}{MeasureTheory.projectiveLimitWithWeakestHypotheses}.

\begin{theorem}[Kolmogorov extension]\label{T1}
  For all $t\in T$, let $E_t$ be a second-countable, complete
  pseudo-extended-metric space and $\mathcal F_t$ the Borel
  $\sigma$-algebra generated by its topology. Let $(P_J)_{J\subseteq_f
    T}$ be a projective family of finite measures and $P$ be
  defined on $\mathcal A := \bigcup_{J \subseteq_f T} \mathcal
  F_J$ given by $P(A) = P_J(\pi_J A)$ for $A\in\mathcal F_J$. Then,
  there is a unique extension of $P$ to $\sigma(\mathcal A)$, which we call projective limit of the family.
\end{theorem}

The uniqueness of the projective limit implies in particular that the finite dimensional distributions of a stochastic process define a unique law.

Note that we extend standard assumptions (see e.g.\ \cite{Billingsley1995,bogachev2007measure}) in two directions. First, we allow that $E$ is a dependent type, i.e.\ for every index $t \in T$, the state space $E_t$ might be different. Second, all $E_t$ are not necessarily second-countable, complete metric spaces (or Polish, i.e.\ second-countable and metrizable through a complete metric), but extended pseudo-metric spaces. Such spaces do not satisfy the frequently used Hausdorff (or T2) property, i.e.\ there can be $x\neq y$ such that all open balls around $x$ and $y$ overlap.
While the first generalization did not require any change in the proof, the second generalization was possible since underlying results in \mathlib were already provided on the same level of generality.
More precisely, the \mathlib lemma \mllink{Topology/UniformSpace/Cauchy}{L637-L638}{isCompact_iff_totallyBounded_isComplete}, which shows that a set $A \subseteq E$ is compact iff it is complete and totally bounded, requires $E$ to be a uniform space.
This is a weaker topological notion than a pseudo-metric space but in the case of a second-countable space (as is the case in our formalization), a uniform space can also be turned into a pseudo-metric space (\mllink{Topology/Metrizable/Uniformity}{L256-L258}{UniformSpace.pseudoMetrizableSpace}).

\subsection{Gaussian Measures and characteristic functions}
\label{S:gaussian}

Our goal is to define a process with Gaussian finite-dimensional distributions.
Prior to our work, \mathlib contained the definition of the one-dimensional normal distribution \mllink{Probability/Distributions/Gaussian/Real}{L198-L201}{ProbabilityTheory.gaussianReal}, but not more general Gaussian measures.
We will denote by $N(m, \sigma^2)$ the normal distribution on $\mathbb{R}$ with mean $m$ and variance $\sigma^2$.
In order to define the law of a Brownian motion on $\mathbb{R}_+$, we need definitions and properties of multivariate Gaussian distributions.
We did not however define Gaussian measures only for that setting, but first introduced a general definition for Banach spaces, and then specialized it to Hilbert and finally finite dimensional spaces.
Our formalization uses many results from \mathlib and we relied on characteristic functions \mllink{MeasureTheory/Measure/CharacteristicFunction}{L322-L325}{MeasureTheory.charFunDual}, as well as the fact that they uniquely determine a probability measure \mllink{MeasureTheory/Measure/CharacteristicFunction}{L435-L439}{MeasureTheory.Measure.ext_of_charFunDual}.

\subsubsection{Gaussian distributions}

We start by recalling the definition of the characteristic function of a general measure. Fix $E$ a real Banach space, and denote by $E'$ its topological dual, i.e. the set of continuous linear forms on $E$. Consider $\mu$ a measure over $E$, and denote by $\mathbf E_\mu$ and $\mathrm{var}_\mu$ the associated expectation and variance operators.

\begin{definition}[\mllink{MeasureTheory/Measure/CharacteristicFunction}{L322-L325}{MeasureTheory.charFunDual}]
	The \emph{dual characteristic function} of $\mu$ is the map $\phi_\mu : E' \to \mathbb{C}$ defined by
	$$\phi_\mu(L) = \int_E e^{i L(x)} \mu(\mathrm dx),$$
	where $i$ is the imaginary unit.
\end{definition}

This function is quite useful because if $E$ is second-countable and equipped with the Borel $\sigma$-algebra, and if $\mu$ is finite, then it is characterized by $\phi_\mu$ (\mllink{MeasureTheory/Measure/CharacteristicFunction}{L435-L439}{MeasureTheory.Measure.ext_of_charFunDual}). Thus it is a great tool to deduce certain properties of the measure and we formalized several results in this sense.
For instance, two random variables are independent if and only if the dual characteristic function of their joint law is the product of the dual characteristic functions of the marginal laws (\mllink{Probability/Independence/CharacteristicFunction}{L58-L64}{ProbabilityTheory.indepFun_iff_charFunDual_prod}).

In the special case where $E$ is a Hilbert space with scalar product $\langle \cdot, \cdot \rangle$, Riesz's representation theorem implies that each $L \in E'$ can be uniquely represented by an element $t \in E$. This motivates the following definition.

\begin{definition}[\mllink{MeasureTheory/Measure/CharacteristicFunction}{L128-L129}{MeasureTheory.charFun}]
  If $E$ is an inner product space, the \emph{characteristic function} of $\mu$ is the map $\psi_\mu : E \to \R$ defined by
  $$\psi_\mu(t) = \int_E e^{i \langle t, x \rangle} \mu(\mathrm dx).$$
\end{definition}

This in some sense is a more convenient way to use the dual characteristic function, as it makes computations easier, and ensures the same characterization property.

Let us now turn to Gaussian measures, which are defined through their one-dimensional projections.

\begin{definition}[\mllink{Probability/Distributions/Gaussian/Basic}{L43-L46}{ProbabilityTheory.IsGaussian}]\label{def:gaussian}
  The measure $\mu$ is called \emph{Gaussian} if for each $L \in E'$, the pushforward measure $L_* \mu$ is the Gaussian measure $N\left(\mathbf E_\mu[L], \mathrm{var}_\mu(L)\right)$.
\end{definition}

In the case of a real Gaussian measure $N(m, \sigma^2)$, the characteristic function is given by
$$\psi_{N(m, \sigma^2)}(t) = \exp(im - \sigma^2/2).$$
This proves in particular that a real Gaussian measure is determined by its mean and covariance. The same statement is actually true for any Gaussian measure, if we generalize the notion of covariance.
As our work primarily focuses on the case where $E$ is a Euclidean space, in what follows we assume that $E$ is a second-countable Hilbert space equipped with its Borel $\sigma$-algebra. However most results remain true if $E$ is a Banach space by only changing definitions to refer to the dual space.

\begin{definition}[\mllink{Probability/Moments/CovarianceBilin}{L48-L52}{ProbabilityTheory.covarianceBilin}]
  The \emph{covariance bilinear form} of a measure $\mu$ with finite second moment $\mathbf{E}_\mu[\Vert \cdot \Vert^2]$ is the bilinear form $C_{\mu} : E \times E \to \R$ defined by
  $$C_\mu(s, t) = \mathbf{cov}_\mu(\langle s, \cdot \rangle, \langle t, \cdot \rangle).$$
\end{definition}

The properties of the covariance allow to easily deduce that $C_\mu$ is actually a continuous and positive semidefinite bilinear form. We have the following characterization of Gaussian measures.

\begin{theorem}[\bmlink{Gaussian/Gaussian}{L207-L109}{ProbabilityTheory.isGaussian_iff_charFun_eq}]\label{thm:gaussian_charFun}
  The measure $\mu$ is Gaussian if and only if its characteristic function is given by
  $$\psi_\mu(t) = \exp(i \langle t, \mathbf E_\mu[\mathrm{id}] \rangle - C_\mu(t, t) / 2).$$
\end{theorem}

This simple description of Gaussian measures allows to easily formalize several important facts about them. For instance, the product $N(0, 1)^{\otimes n}$ of $n \in \mathbb{N}$ times many Gaussian measures $N(0, 1)$ is a Gaussian measure.
Another result is that if two random variables which are jointly Gaussian have covariance zero, then they are independent.

Theorem~\ref{thm:gaussian_charFun} used implicitly the fact that Gaussian measures have finite first and second moments $\mathbf E_\mu[\lVert \cdot \rVert]$ and $\mathbf{E}_\mu[\lVert \cdot \rVert^2]$.
While we would only need to prove this result in finite dimension to define a Brownian motion on $\mathbb{R}_+$, we proved a more general theorem applicable to Banach spaces, called Fernique's theorem \cite{fernique1970integrabilite} (\mllink{Probability/Distributions/Gaussian/Fernique}{L161-L164}{ProbabilityTheory.IsGaussian.exists_integrable_exp_sq}).
It states that if $\mu$ is a Gaussian measure over a Banach space $E$, then there exists $C > 0$ such that
\begin{align*}
  \int_E \exp\left( C \|x\|^2 \right) \mu(\mathrm dx) < \infty.
\end{align*}
This implies that all moments of a Gaussian measure in a second-countable Banach space are finite.

\paragraph{On Euclidean spaces}

As stated before, any Gaussian measure is uniquely determined by its mean and covariance. A follow-up question is to know whether, given $m \in E$ and $C$ a continuous and positive semidefinite bilinear form over $E$, there exists a Gaussian measure with mean $m$ and covariance $C$.
The answer is yes if $E$ is a Euclidean space, and we formalized this construction in the project.
As a bilinear form over a finite-dimensional vector-space can be represented by a matrix, we will use the \emph{covariance matrix} to represent the covariance of a measure.

Assume now that $E$ is a Euclidean space with a basis $(e_1, \dots, e_n)$. We define the \emph{standard Gaussian measure} (\bmlink{Gaussian/MultivariateGaussian}{L33-L37}{ProbabilityTheory.stdGaussian}) as the pushforward measure
$$\left(x \mapsto \sum_{i=1}^n x_i e_i\right)_* N(0, 1)^{\otimes n}.$$

As stated above, $N(0, 1)^{\otimes n}$ is a Gaussian measure because it is a product of Gaussian measures, and the map $x \mapsto \sum_{i=1}^n x_i e_i$ is a linear form over $E$, so this measure is indeed Gaussian. It is easy to see that its mean is $0$ and its covariance matrix is the identity matrix $I_n$.
Also, the measure we obtain does not depend on the choice of the basis.
We denote it $N(0, I_n)$.

Now from Definition~\ref{def:gaussian} it is clear that the pushforward of a Gaussian measure by a continuous affine map is again a Gaussian measure. In particular, given $A$ an $n \times n$ matrix, $b \in E$ and $\mu$ a Gaussian measure over $E$ with mean $m$ and covariance $C$, we have that $(x \mapsto Ax + b)_* \mu$ is a Gaussian measure with mean $Am + b$ and covariance matrix $A^\top C A$.
Consider thus $b \in E$ and $A$ a positive semidefinite $n \times n$ matrix.
There exist $S \in \mathcal{M}_n(\R)$ such that $A = S^\top S$, and thus $(x \mapsto Sx + b)_* N(0, I_n)$ is a Gaussian measure with mean $m$ and covariance matrix $A$. This proves that for any $m \in E$ and any covariance matrix $A$, there exist a Gaussian measure over $E$ with mean $m$ and covariance matrix $A$. This measure is denoted $N(m, A)$, and is formalized as \bmlink{Gaussian/MultivariateGaussian}{L150-L155}{ProbabilityTheory.multivariateGaussian}.

\subsubsection{Projective family of finite-dimensional Gaussian measures}

The law of the standard Brownian motion on $\mathbb{R}_+$ is a measure on $\mathbb{R}^{\mathbb{R}_+}$, whose finite-dimensional distributions are specific multivariate Gaussian measures.
We build it as the projective limit of these finite-dimensional distributions, using the Kolmogorov extension theorem (Theorem~\ref{T1}).
In order to apply the theorem, we need to show that these finite-dimensional distributions form a projective family of measures.

Let then $P_J$ be the measure over $\mathbb{R}^J$ defined as the multivariate Gaussian measure $N(0, C_J)$, where $C_J$ is the covariance matrix defined by $(C_J)_{ij} = t_i \wedge t_j$ for $J = \{t_1, \dots, t_n\} \subseteq \mathbb{R}_+$.
To show that this Gaussian measure exists we need to prove that $C_J$ is positive semidefinite.
For this, we rely on Gram matrices, which are based on inner product spaces. In such a space $E$, for $v_i,...,v_n \in E$, define a matrix $G$ by $G_{ij} := \langle v_i, v_j\rangle$.
Then, for $t :=(t_1,...,t_n) \in \R^n$,
\begin{align*}
  t^\top G t = \sum_{i,j} t_i \langle v_i, v_j\rangle t_j = \Big\langle \sum_i t_i v_i, \sum_j t_j v_j\Big\rangle = \left\Vert \sum_i t_i v_i \right\Vert^2 \geq 0
  \: ,
\end{align*}
and thus the Gram matrix $G$ is positive semidefinite.
Let $v_i = 1_{[0,t_i]}$ in the space of $L^2$ integrable functions with respect to Lebesgue integral. Then,
$$ \langle v_i, v_j \rangle = \int 1_{[0,t_i]}(x) 1_{[0,t_j]}(x) dx = \int 1_{[0,t_i \wedge t_j]}(x) dx = t_i \wedge t_j.$$
So, $C_J$ from above is a Gram matrix, which is positive semidefinite.

In order to prove that $(P_J)_{J \subseteq_f \mathbb{R}_+}$ is projective, we need to show that the pushforward of a Gaussian measure $P_J = N(0, C_J)$ by the projection $\pi_H^J : \mathbb{R}^J \to \mathbb{R}^H$ is equal to $P_H = N(0, C_H)$ for all $H \subseteq J \subseteq_f \mathbb{R}_+$.
Since we have proved that the pushforward of a Gaussian measure by a continuous linear map is again a Gaussian measure, it suffices to compute its mean and covariance, which finishes the proof of projectivity.

We have thus constructed a measure on $\mathbb{R}^{\mathbb{R}_+}$.
We can define the canonical process on $\mathbb{R}^{\mathbb{R}_+}$, which is the process $(X_t)_{t \in \mathbb{R}_+}$ with $X_t(\omega) = \omega(t)$ for all $\omega \in \mathbb{R}^{\mathbb{R}_+}$.
By construction, this process has the right finite-dimensional distributions to be a Brownian motion.
However, its paths are not necessarily continuous.

\subsection{The Kolmogorov-Chentsov Theorem}
\label{S:continuity}

The Kolmogorov-Chentsov theorem gives the existence of a Hölder continuous modification of a stochastic process, provided that the process satisfies a certain moment condition called the Kolmogorov condition.

\subsubsection{The Kolmogorov condition}

Suppose that the index set $T$ and the value space $E$ of a stochastic process are extended pseudo-metric spaces.
That is, $T$ is equipped with a ``distance'' $d_T$ with values in $[0,+\infty]$ (infinity included), which is reflexive, symmetric and satisfies the triangle inequality, and a similar distance $d_E$ is defined on $E$. ``extended'' refers to the fact that the distance can take the value $+\infty$, and ``pseudo'' to the possibility of non-identical points having distance zero.
$T$ will denote the index set of a stochastic process, and $E$ the space of values of that process.
The space $E$ is equipped with the Borel $\sigma$-algebra $\mathcal{B}(E)$ generated by the open sets of $E$.

The key assumption of the Kolmogorov-Chentsov theorem that gives the existence of a continuous modification is that the process satisfies the Kolmogorov condition, \mllink{Probability/Process/Kolmogorov}{L52-L60}{ProbabilityTheory.IsKolmogorovProcess}.

\begin{definition}[Kolmogorov condition]\label{def:kolmogorov_condition}
Let $p, q, M$ be non-negative real numbers with $p,q>0$.
A stochastic process $(X_t)_{t \in T}$ is said to satisfy the Kolmogorov condition for exponents $(p, q)$ with constant $M$ if for all $s, t \in T$, the pair $(X_t, X_t) : \Omega \to E \times E$ is measurable for the Borel $\sigma$-algebra on $E \times E$, and
\begin{align*}
  \mathbf{E}[d_E(X_s, X_t)^p] \le M d_T(s, t)^q
  \: .
\end{align*}
\end{definition}

The measurability condition ensures that the distance $d_E(X_s, X_t)$ is a measurable.
If $E$ was second-countable, we could simply assume that $X_s$ is measurable for all $s \in T$, and the measurability of the pair would follow.
This is because for a second-countable space $E$, $\mathcal{B}(E \times E)$ is the product of the Borel $\sigma$-algebras on $E$, which is not true in general.

We remark that the inequality in Definition~\ref{def:kolmogorov_condition} remains true for any modification of the process $X$, and that Theorem~\ref{thm:kolmogorov_chentsov} also holds for a modification of a process satisfying the Kolmogorov condition.
Therefore, we introduce a definition \lstinline|IsAEKolmogorovProcess| for the property of being a modification of a process satisfying the Kolmogorov condition, and prove the theorem under that hypothesis.

\subsubsection{The main theorem}

We formalize a recent and general version of the Kolmogorov-Chentsov theorem \cite[Theorem 1]{kratschmer2023kolmogorov}, which applies to extended pseudo-metric spaces under a covering assumption on the index set $T$.
We say that $T$ has bounded covering number with constant $c > 0$ and exponent $d \ge 0$ if for all $\varepsilon \in (0, \mathrm{diam}(T)]$, $T$ can be covered by at most $c \varepsilon^{-d}$ balls of radius $\varepsilon$ (see Definition~\ref{def:bounded_covering_number} below). Note that such a space is totally bounded by construction.

\begin{theorem}[Kolmogorov-Chentsov]\label{thm:kolmogorov_chentsov}
Suppose that the index set $T$ has bounded covering number with constant $c>0$ and exponent $d > 0$.
Let $(X_t)_{t \in T}$ be a stochastic process that satisfies the Kolmogorov condition for exponents $(p,q)$ with constant $M$, with $q > d$ and $p > 0$.
Then for all $\beta \in(0, (q - d)/p)$ there exists a finite constant $L(T, c, d, p, q, \beta)$ such that for every countable subset $T' \subseteq T$,
\begin{align*}
  \mathbf{E}\left[ \sup_{s, t \in T'} \frac{d_E(X_s, X_t)^p}{d_T(s, t)^{\beta p}} \right]
  \le M L(T, c, d, p, q, \beta)
  \: .
\end{align*}
If furthermore $E$ is complete and $T$ is second-countable, then the process $X$ has a modification with Hölder continuous paths of exponent $\beta$ for all $\beta \in (0, (q - d)/p)$.
\end{theorem}

In our formalization, the two parts of this theorem correspond to the statements \bmlink{Continuity/KolmogorovChentsovInequality}{L405-L410}{ProbabilityTheory.countable_kolmogorov_chentsov} for the inequality on the expectation and \bmlink{Continuity/KolmogorovChentsov}{L1216-L1221}{ProbabilityTheory.exists_modification_holder} for the existence of a Hölder modification.
We adopt here \Lean's convention that $0/0 = 0$, which means that the expression under the expectation is well-defined even if we allow $d_T(s, t) = 0$ in the supremum.
Indeed, if $d_T(s, t) = 0$, then the Kolmogorov condition implies that $d_E(X_s, X_t) = 0$ almost surely, so the ratio is then $0/0 = 0$ and the expectation is well-defined.

Note that $T$ is said to be an \emph{extended} pseudo-metric space, but the bounded covering number condition implies that $T$ has finite diameter, hence is a pseudo-metric space.
The \Lean code is nonetheless written for a distance with values in $[0,+\infty]$, for two reasons: first the distance on $E$ has that type and using the same types for both is easier, and second it allows us to work with \mathlib's Lebesgue integral (which expect that co-domain for the functions we integrate).
Working with the Lebesgue integral is easier than using \mathlib's Bochner integral since manipulating the latter often requires proving integrability of the functions we consider.
Therefore, from the beginning of the project we made the design choice to consider any non-negative function as taking values in $[0,+\infty]$, and to use the Lebesgue integral whenever possible.

Our reference \cite{kratschmer2023kolmogorov} proves the theorem for $T$ and $E$ metric spaces.
The most notable difference is that we allow pseudo-metrics, for which $d_T(s, t)$ can be zero for $s \ne t$.
The proof we implement follows that paper, with minor changes.
The first one is that the authors of \cite{kratschmer2023kolmogorov} used the Minkowski inequality in a part of the proof, which requires $p \ge 1$, and did not mention how to handle the (easier) case $p < 1$, although the theorem is also valid there.
Beyond fixing this minor omission by adding a separate argument for $p \in (0,1)$, the second change we make is to ensure that all exponents we use are natural numbers, in order to simplify the formalization.
Whenever the proof of \cite{kratschmer2023kolmogorov} uses a geometric grid $2^{-n_1}, 2^{-n_1-1},\ldots, 2^{-n_2}$ with exponents $n_1 < n_2 \in \mathbb{Z}$ possibly negative, we introduce a multiplicative factor $\varepsilon_0$ and use a grid of the form $\varepsilon_0 2^{-m_1}, \ldots, \varepsilon_0 2^{-m_2}$ with $m_1, m_2 \in \mathbb{N}$.
\mathlib has more support for sums over natural numbers, and natural number exponents are easier to manipulate as for example the monotonicity of $x \mapsto x^r$ for $r \in \mathbb{Z}$ depends on the sign of $r$ but is always the same for $r \in \mathbb{N}$. In practice, using natural number exponents means that we don't have to prove as many side conditions in the formalization.
Those changes do not impact the main ideas of the proof, but the constants we report are different from those of \cite{kratschmer2023kolmogorov}.



\paragraph{Localized theorem}

In order to obtain a Brownian motion, we will want to apply Theorem~\ref{thm:kolmogorov_chentsov} to a process defined on $T = \mathbb{R}_+$, which is not bounded and does not have bounded covering number.
We thus need an extension of the theorem.
The following localized version of the theorem fulfills this need, and also uses several sequences of exponents $(p_m, q_m)$ instead of a single pair $(p, q)$ to improve the exponents in the local Hölder continuity of the resulting process.
The conclusion of the theorem is not Hölder continuity of the paths, but local Hölder continuity.
A function $f : T \to E$ is said to be locally Hölder continuous of order $\gamma$ if for all $t \in T$, there exists a neighborhood $U$ of $t$ such that $f$ is $\gamma$-Hölder continuous on $U$.

\begin{theorem}[\bmlink{Continuity/KolmogorovChentsov}{L1362-L1368}{ProbabilityTheory.exists_modification_holder_iSup}]\label{thm:localized_holder_modification_sup}
Suppose that $T$ can be covered by an increasing sequence of totally bounded open subsets $(T_n)_{n \in \mathbb{N}}$ such that each $T_n$ has bounded covering number with constant $c_n > 0$ and exponent $d > 0$ (the same exponent for all $n$).
Let $(p_m, q_m)_{m \in \mathbb{N}}$ be a sequence of pairs of positive numbers such that $q_m > d$ for all $m \in \mathbb{N}$.
Let $(X_t)_{t \in T}$ be a process that satisfies the Kolmogorov condition with exponents $(p_m, q_m)$ for all $m \in \mathbb{N}$.
Then $X$ has a modification $Y$ such that the paths of $Y$ are locally Hölder continuous of all orders $\gamma \in \left(0, \sup_m \frac{q_m - d}{p_m}\right)$.
\end{theorem}

The hypothesis on $T$ is satisfied for $T = \mathbb{R}_+$, with $T_n = [0,n)$ (which is open in $\mathbb{R}_+$), for $d = 1$.
Although we use it only for $T = \mathbb{R}_+$, Theorem~\ref{thm:localized_holder_modification_sup} can be used in higher dimension and in more involved cases: the authors of \cite{kratschmer2023kolmogorov} show how to apply it to subsets of $m$-dimensional Riemannian manifolds.

We don't describe the proof of Theorem~\ref{thm:localized_holder_modification_sup} in detail here (we refer the reader to \cite{kratschmer2023kolmogorov}, or to the code): it consists in applying the same type of arguments as in the proof of the second part of Theorem~\ref{thm:kolmogorov_chentsov} to build separately modifications for each set $T_n$ and exponents ($p_m, q_m)$ and then combine them into one process (see section~\ref{sub:holder_process}).

\subsubsection{Covers and chaining}

The proof of Theorem~\ref{thm:kolmogorov_chentsov} uses discretizations of the space $T$ at different scales, using covers, and uses a chaining argument to relate the supremum of distances over $T$ to suprema over those discretizations.
Covers and chaining are important tools in the study of suprema of stochastic processes, in particular in bounds on their expectation or on their tail probabilities \cite{talagrand2022upper, vershynin2018high}.
Part of the reason for choosing to formalize the Kolmogorov-Chentsov theorem as presented in \cite{kratschmer2023kolmogorov} is that it uses these techniques, which can be reused in other contexts.

\paragraph{Covers}

An (internal) $\varepsilon$-cover of $T$ is a subset $S \subseteq T$ such that for all $t \in T$, there exists $s \in S$ such that $d_T(s, t) \le \varepsilon$.
If $T$ has bounded diameter we can find an $\varepsilon$-cover which is finite, and the covering number of $T$ for $\varepsilon > 0$ is the minimal cardinality of an $\varepsilon$-cover of $T$, denoted by $N_\varepsilon(T)$.
We call a cover with that cardinality a minimal $\varepsilon$-cover.
Our formalization also includes definitions of external covers and covering numbers, packing numbers and inequalities between these numbers and the volume of balls. We don't detail those here, but the relations between those quantities are helpful to prove properties of covering numbers.

The main hypothesis of the Kolmogorov-Chentsov theorem on the index set is that $T$ has bounded covering number, defined as follows.

\begin{definition}[\bmlink{Continuity/HasBoundedInternalCoveringNumber}{L18-L19}{HasBoundedInternalCoveringNumber}]\label{def:bounded_covering_number}
A set $T$ has bounded covering number with constant $c > 0$ and exponent $d \ge 0$ if for all $\varepsilon \in (0, \mathrm{diam}(T)]$, $N_\varepsilon(T) \le c \varepsilon^{-d}$.
\end{definition}

The exponent is a measure of dimensionality of the space $T$.
For instance, for $T = [0,1]^m$ with the Euclidean distance, $T$ has bounded covering number with exponent $d = m$.
Covering numbers are not monotone for set inclusion, but if $A \subseteq B$ then $N_\varepsilon(A) \le N_{\varepsilon/2}(B)$, which implies the following lemma on bounded covering numbers of subsets.

\begin{lemma}[\bmlink{Continuity/HasBoundedInternalCoveringNumber}{L48-L50}{HasBoundedInternalCoveringNumber.subset}]\label{lem:covering_number_subset}
Let $T$ be a set with bounded covering number with constant $c$ and exponent $d$, and let $T' \subseteq T$ be a subset of $T$.
Then $T'$ has bounded covering number with the constant $2^d c$ and exponent $d$.
\end{lemma}

\paragraph{Chaining}

We will consider minimal covers $(C_n)_{n \in \mathbb{N}}$ for a sequence $(\varepsilon_n)_{n \in \mathbb{N}}$ of positive real numbers decreasing to $0$, which in the proof will be chosen as $\varepsilon_n = \varepsilon_0 2^{-n}$.
The proof of Theorem~\ref{thm:kolmogorov_chentsov} will involve bounding a supremum $\mathbf{E}\left[\sup_{s, t \in C_N, \: d_T(s, t) \le \delta}d_E(X_s, X_t)^p\right]$ for a fine scale cover $C_N$.
In order to do so, we will move to a coarser scale $m < N$.
For $x \in C_N$, we define a chaining sequence $\bar{x}_0, \bar{x}_1, \ldots, \bar{x}_N$ by letting $\bar{x}_N = x$ and for $n < N$, $\bar{x}_n$ be a point in $C_n$ closest to $\bar{x}_{n+1}$.
That sequence satisfies $d_T(\bar{x}_n, \bar{x}_{n+1}) \le \varepsilon_n$ for all $n < N$, by the cover property of $C_n$.

The change of scale from $C_N$ to $C_m$ is done through the following decomposition, due to the triangle inequality and the fact that $(a + b)^p \le 2^p (a^p + b^p)$ for $a, b \ge 0$ (\bmlink{Continuity/Chaining}{L257-L262}{scale_change_rpow}):
\begin{align*}
  \sup_{s, t \in C_N, d_T(s, t) \le \delta} d_E(X_s, X_t)^p
  &\le 2^p \sup_{s, t \in C_N, d_T(s, t) \le \delta} d_E(X_{\bar{s}_m}, X_{\bar{t}_m})^p + 4^p \sup_{s \in C_N} d_E(X_s, X_{\bar{s}_m})^p
  \: .
\end{align*}
The first term on the right hand side is a supremum with points in the coarser cover $C_m$, while the second term quantifies the distance along the chaining sequence.
The first term will be bounded by $2^p \sup_{s, t \in C_m, d_T(s, t) \le \delta'} d_E(X_s, X_t)^p$ for a constant $\delta'$ related to $\delta$ and the $\varepsilon_n$'s.
It thus has the same form as the original supremum, but on a coarser scale.
By the triangle inequality and the cover property, for $x, y \in C_N$ with $d_T(x, y) \le \delta$, $d_T(\bar{x}_m, \bar{y}_m) \le \delta + 2 \sum_{i=m}^{N-1}\varepsilon_i$.
For a geometric grid $\varepsilon_n = \varepsilon_0 2^{-n}$ and $m$ such that $\varepsilon_m \le \delta \le 4 \varepsilon_m = 4 \varepsilon_0 2^{-m}$, we obtain $d_T(\bar{x}_m, \bar{y}_m) \le \varepsilon_0 2^{-m+3}$, such that
\begin{align}
  \sup_{s, t \in C_N, d_T(s, t) \le \delta} d_E(X_s, X_t)^p
  &\le 2^p \sup_{s', t' \in C_m, d_T(s', t') \le \varepsilon_0 2^{-m+3}} d_E(X_{s'}, X_{t'})^p + 4^p \sup_{s \in C_N} d_E(X_s, X_{\bar{s}_m})^p
  \: . \label{eq:scale_change}
\end{align}

\paragraph{Pair reduction lemma}

In the proof of Theorem~\ref{thm:kolmogorov_chentsov} we will need to control the expectation of the first term of the right hand side of \eqref{eq:scale_change}.
Essentially, we bound the expectation of a supremum over pairs by the number of pairs multiplied by an upper bound for the expectation of each distance.
A bound on the cardinality of $\{(s, t) \in C_m^2 \mid d_T(s, t) \le \delta'\}$ would be $\lvert C_m \rvert^2$, which as we explain in the main proof below is too large for our needs.
The pair reduction lemma below allows us to reduce the number of pairs we need to consider to a quantity of order $\lvert C_m \rvert$ instead of quadratic, at the cost of increasing the distance threshold from $\delta'$ to a larger value $\delta' \log \lvert C_m \rvert$.

\begin{lemma}[\cite{kratschmer2023kolmogorov} Lemma 6.1, extension of \cite{talagrand2014}, B.2.7, \mllink{Topology/EMetricSpace/PairReduction}{L468-L482}{EMetric.pair_reduction}]\label{lem:pair_reduction}
Let $T$ be a pseudo-metric space, $c \ge 0$, $J$ a finite subset of $T$, $a \ge 0$ and $n \in \mathbb{N}$ such that $\lvert J \rvert \le a^n$.
Then there exists a set $K \subseteq J^2$ such that for any pseudo-metric space $E$ and any map $f: T \to E$~,
\begin{enumerate}
  \item $\lvert K \rvert \le a \lvert J \rvert$~,
  \item for all $(s, t) \in J^2$, $d_T(s, t) \le cn$~,
  \item $\sup_{s, t \in J, \: d_T(s, t) \le c} d_E(f(s), f(t)) \le 2 \sup_{(s, t) \in K} d_E(f(s), f(t))$~.
\end{enumerate}
\end{lemma}


\subsubsection{Proof of the main inequality}

Under the assumptions of Theorem~\ref{thm:kolmogorov_chentsov}, we prove that for all countable subset $T' \subseteq T$,
\begin{align*}
  \mathbf{E}\left[ \sup_{s, t \in T'} \frac{d_E(X_s, X_t)^p}{d_T(s, t)^{\beta p}} \right]
  \le M L(T, c, d, p, q, \beta)
  \: ,
\end{align*}
for a finite constant $L(T, c, d, p, q, \beta)$ to be determined.
This section does not cover all details of the proof, for which we refer to \cite{kratschmer2023kolmogorov} and to our blueprint. We only highlight the main steps, key arguments, and places that are relevant to the formalization specifically.
In particular, we explain why the pair reduction Lemma~\ref{lem:pair_reduction} is needed in the proof.

First, it suffices to show the inequality for finite subsets $J \subseteq T$, since the countable case follows by monotone convergence.
Since $J \subseteq T$, it also has bounded covering number with the same exponent $d$ and a constant $2^d c$ (Lemma~\ref{lem:covering_number_subset}).

The next step is to discretize the possible distances between points in $J$.
For $k \in \mathbb{N}$, let $\eta_k = 2^{-k}(\mathrm{diam}(T) + 1) \ge 2^{-k}$.
We write the set of pairs of points in $J$ as the union of the sets $\{(s, t) \in J^2 \mid d_T(s, t) \in (\eta_k, 2\eta_k]\}$.
Then we can show
\begin{align*}
  \mathbf{E}\left[ \sup_{s, t \in J} \frac{d_E(X_s, X_t)^p}{d_T(s, t)^{\beta p}} \right]
  &\le \sum_{k=0}^\infty 2^{k \beta p} \mathbf{E}\left[ \sup_{\substack{s, t \in J \\ d_T(s, t) \le 2\eta_k}} d_E(X_s, X_t)^p \right]
  \: .
\end{align*}
We thus want to find an upper bound, for $\delta > 0$, on
$
  \mathbf{E}\left[ \sup_{s, t \in J,\: d_T(s, t) \le \delta} d_E(X_s, X_t)^p \right]
  \: .
$
We will then check that the series we obtain when summing over $k$ converges when $\beta < (q - d)/p$, and its value will be the constant $L(T, c, d, p, q, \beta)$.

\begin{lemma}[\bmlink{Continuity/IsKolmogorovProcess}{L994-L1000}{ProbabilityTheory.finite_set_bound_of_edist_le}]\label{lem:finite_set_bound_of_dist_le}
For $J$ a finite set with bounded covering number with constant $c > 0$ and exponent $d > 0$, for $X$ a stochastic process that satisfies the Kolmogorov condition for exponents $(p, q)$ with constant $M$, with $q > d$ and $p > 0$, and for $\delta > 0$, we have
\begin{align*}
  \mathbf{E}\left[ \sup_{\substack{s, t \in J \\ d_T(s, t) \le \delta}} d_E(X_s, X_t)^p \right]
  &\le 2^{2p+4q+1} M c \delta^{q-d}\left( 4^d \left(\max\left\{0, \log_2 \left(c 4^d \delta^{-d}\right)\right\}\right)^q + R_p \right)
  \: ,
  \\
  \text{ with } R_p
  &= \max\left\{\frac{1}{2^{q-d} - 1}, \frac{1}{(2^{(q-d)/p} - 1)^p}\right\}
  \: .
\end{align*}
\end{lemma}

We apply that lemma for $\delta$ of order $2^{-k}$ and we obtain for the constant $L$ a sum over $k$ of $2^{k(\beta p - (q - d))}$ multiplied by an expression that grows like $k^q$.
That series converges if and only if $\beta < (q - d)/p$ (a fact that was surprisingly tedious to prove using current tactics from \mathlib, although recent developments on a tactic for computing limits may help in the future).
To obtain that exact condition on $\beta$, it is thus important to get the factor $\delta^{q-d}$ in the lemma above, and not a worse exponent.
This is where the chaining argument is used, as seen further down in the proof.

The remainder of this section is devoted to the proof of Lemma~\ref{lem:finite_set_bound_of_dist_le}.
We note first that we can restrict the supremum in the expectation to pairs $(s, t) \in J^2$ such that $d_T(s, t) > 0$, since if $d_T(s, t) = 0$, then by the Kolmogorov condition, $d_E(X_s, X_t) = 0$ almost surely.
Let $\varepsilon_0 = \mathrm{diam}(T)$ and for $n \in \mathbb{N}$, let $\varepsilon_n = 2^{-n} \varepsilon_0$.
Since $J$ is finite, there exists $N \in \mathbb{N}$ such that $\varepsilon_N < \min_{s, t \in J, \: d_T(s, t) > 0} d_T(s, t)$. Let $N$ be the smallest such integer.
We can assume that $\delta \ge \varepsilon_N$, since otherwise the supremum is zero.
The supremum over the set $J$ is thus equal to the supremum over minimal cover $C_N$ of $J$ by balls of radius $\varepsilon_N$.
Let more generally $C_n$ be a minimal cover of $J$ by balls of radius $\varepsilon_n$ for all $n \in \{0, \dots, N\}$.

We change scale using Equation~\eqref{eq:scale_change}, with $m \in \{0, \dots, N-1\}$ chosen such that $\varepsilon_m \le \delta \le 4 \varepsilon_m$.
Note that this needs $\delta \le 4 \mathrm{diam}(T)$, and we omit here the case $\delta > 4 \mathrm{diam}(T)$ for which we change the scale to the singleton cover $C_0$ and where the proof is a simpler version of what follows (\bmlink{Continuity/IsKolmogorovProcess}{L674-L678}{ProbabilityTheory.finite_set_bound_of_edist_le_of_diam_le}).
\begin{align*}
  \mathbf{E}\left[\sup_{s, t \in C_N, d_T(s, t) \le \delta} d_E(X_s, X_t)^p\right]
  &\le 2^p \mathbf{E}\left[\sup_{s', t' \in C_m, d_T(s', t') \le 8 \varepsilon_m} d_E(X_{s'}, X_{t'})^p\right] + 4^p \mathbf{E}\left[\sup_{s \in C_N} d_E(X_s, X_{\bar{s}_m})^p\right]
  \: .
\end{align*}
We now bound separately the two terms in that sum, which correspond respectively to distances between points in the coarser cover $C_m$, and to distances along the chaining sequences.

\textbf{First term.}
The main consequence of the Kolmogorov condition that we will use in the proof is that for $K$ a finite set of pairs of points in $T$ such that for all $(s,t) \in K$, $d_T(s, t) \le \varepsilon$, we have
\begin{align}
  \mathbf{E}\left[ \sup_{(s,t) \in K} d_E(X_s, X_t)^p \right]
  \le \sum_{(s,t) \in K} \mathbf{E}\left[ d_E(X_s, X_t)^p \right]
  \le M \lvert K \rvert \varepsilon^q
  \: .\label{eq:kolmogorov_condition_finite_set}
\end{align}

Let us first present a naive upper bound on the expectation of the supremum, which will highlight the need for more advanced techniques.
We apply the inequality~\eqref{eq:kolmogorov_condition_finite_set} for $K = \{(s,t) \in C_m^2 \mid d_T(s, t) \le 8 \varepsilon_m\}$ with cardinal $\lvert K \rvert \le \lvert C_m \rvert^2$, and then the bounded covering number assumption:
\begin{align*}
  \mathbf{E}\left[ \sup_{\substack{s, t \in C_m \\ d_T(s, t) \le 8 \varepsilon_m}} d_E(X_s, X_t)^p \right]
  &\le M \lvert C_m \rvert^2 (8 \varepsilon_m)^q
  \le M c^2 \varepsilon_m^{-2d} (8 \varepsilon_m)^q
  \: .
\end{align*}
Since $\varepsilon_m \le \delta \le 4 \varepsilon_m$, we obtain a bound of order $\delta^{q - 2d}$, which is not sufficient to conclude the proof with the right exponent $\delta^{q - d}$.
The issue is the square on the cardinal of the cover $C_m$, and this is where the pair reduction lemma (Lemma~\ref{lem:pair_reduction}) comes into play.
With that lemma applied to $a = 2$, $n = \lceil \log_2 J \rceil$ and $c = 8 \varepsilon_m$, we obtain a set $K \subseteq C_m^2$ with $\lvert K \rvert \le 2 \lvert C_m \rvert$ and $\sup_{(s, t) \in K} d_T(s, t) \le 8 n \varepsilon_m$, such that an application of \eqref{eq:kolmogorov_condition_finite_set} gives
\begin{align*}
  \mathbf{E}\left[ \sup_{\substack{s, t \in C_m \\ d_T(s, t) \le 8 \varepsilon_m}} d_E(X_s, X_t)^p \right]
  &\le 2^p \mathbf{E}\left[ \sup_{s, t \in K} d_E(X_s, X_t)^p \right]
  \le 2^p M \lvert K \rvert (8 n \varepsilon_m)^q
  \: .
\end{align*}
Omitting constants, that is now of order
$\lvert C_m \rvert \varepsilon_m^q (\log_2 \lvert C_m\rvert)^q \approx \delta^{q - d} \left(\log_2 \delta^{-d} \right)^q
$, which has the right exponent on $\delta$ to conclude the proof.

\textbf{Second term.}

To obtain an upper bound on $\mathbf{E}\left[\sup_{s \in C_N} d_E(X_s, X_{\bar{s}_m})^p\right]$, we consider separately the cases $p \in (0, 1]$ and $p \ge 1$.
For $p \le 1$ we use that the power function is sub-additive: for $a, b \ge 0$, $(a + b)^p \le a^p + b^p$.
\begin{align*}
  \mathbf{E}\left[\sup_{s \in C_N} d_E(X_s, X_{\bar{s}_m})^p\right]
  &\le \mathbf{E}\left[\sup_{s \in C_N} \sum_{i=m}^{N-1} d_E(X_{\bar{s}_{i+1}}, X_{\bar{s}_i})^p\right]
  \le \sum_{i=m}^{N-1} \mathbf{E}\left[\sup_{s \in C_N} d_E(X_{\bar{s}_{i+1}}, X_{\bar{s}_i})^p\right]
  \: .
\end{align*}
It remains to control the expectation inside the sum, which we do using the Kolmogorov condition though Equation~\eqref{eq:kolmogorov_condition_finite_set} with $K = \{(\bar{s}_{i+1}, \bar{s}_i) \mid s \in C_N\}$. That set has cardinal at most $\lvert C_{i+1} \rvert$ (the number of possible values for $\bar{s}_{i+1}$, since $\bar{s}_i$ is a function of $\bar{s}_{i+1}$).
\begin{align*}
  \mathbf{E}\left[\sup_{s \in C_N} d_E(X_{\bar{s}_{i+1}}, X_{\bar{s}_i})^p\right]
  &\le M \lvert C_{i+1} \rvert \varepsilon_i^q
  \le M c \varepsilon_{i+1}^{-d} \varepsilon_i^{q}
  \: .
\end{align*}
We omit the computations, but the sum of those expressions is bounded by the $R_p$ term in the statement of the lemma.
For $p \ge 1$, we use the Minkowski inequality (see \cite{kratschmer2023kolmogorov} for details) and obtain
\begin{align*}
  \mathbf{E}\left[\sup_{s \in C_N} d_E(X_s, X_{\bar{s}_m})^p\right]
  &\le \left( \sum_{i=m}^{N-1} \left( \mathbf{E}\left[\sup_{s \in C_N} d_E(X_{\bar{s}_{i+1}}, X_{\bar{s}_i})^p\right] \right)^{1/p} \right)^p
  \: .
\end{align*}
The expectation inside the sum is bounded as before, and again we get something bounded by the $R_p$ term in Lemma~\ref{lem:finite_set_bound_of_dist_le}.
Only the case $p \ge 1$ was treated in \cite{kratschmer2023kolmogorov}, although the final result is presented for all $p > 0$.
Accordingly, they write their lemma with only the part $(2^{(q-d)/p} - 1)^{-p}$ in the definition of $R_p$, which we correct to include the $p \le 1$ case.
We finally note that keeping track of how the constants evolve through the proof is straightforward in \Lean, in contrast to for example the latex proof which we wrote in the blueprint, in which the constants needed to be modified repeatedly due to small mistakes.
Also any improvement in a constant in a \Lean proof is easy to propagate through the rest of the proof, since the system automatically raises errors where the new value invalidates a computation.

\subsubsection{Building a Hölder continuous process}
\label{sub:holder_process}

We discuss the proof of the second part of Theorem~\ref{thm:kolmogorov_chentsov}, which states that for $E$ complete and $T$ second-countable, the process $X$ has a modification with Hölder continuous paths.
We first prove the result for a fixed $\beta \in (0, (q - d)/p)$, and will then deduce it for all $\beta$ in that interval simultaneously.
An advantage of the formal proof is that it is easy to inspect the code to see which assumptions are used where.
This part of the theorem requires $T$ to be second-countable because we first take a countable dense subset $T'$ of $T$, on which we apply the main inequality of Theorem~\ref{thm:kolmogorov_chentsov}.
We obtain that an expectation is finite, hence that the quantity inside the expectation is almost surely finite:
\begin{align*}
  \sup_{s, t \in T'} \frac{d_E(X_s, X_t)^p}{d_T(s, t)^{\beta p}} < +\infty
  \quad \text{almost surely.}
\end{align*}

On the event $A$ that the supremum is finite, $(X_t)_{t \in T'}$ has Hölder continuous paths of order $\beta$.
Let $x_0 \in E$ be arbitrary and let $Y$ be the process defined by
\begin{align*}
  Y_t(\omega)
  &= \begin{cases}
    \lim_{s \to t, s \in T'} X_s(\omega) & \text{if } \omega \in A \: , \\
    x_0 & \text{otherwise .}
  \end{cases}
\end{align*}
The limits in that definition are the reason for the completeness assumption on $E$.
Then $Y$ has Hölder continuous paths of order $\beta$ and is measurable.
We finally need to show that $Y$ is a modification of $X$.
It is actually only a ``distance zero modification'', by which we mean that $d_E(X_t, Y_t) = 0$ almost surely.
Since $E$ is only a pseudo-metric space, $Y$ might not satisfy the modification property $Y_t = X_t$ almost surely.
We refer to \cite{kratschmer2023kolmogorov} for the proof that $Y$ is a distance zero modification of $X$\footnote{They use a metric space assumption for $E$ and prove that $Y$ is a true modification: replacing equalities with distance zero statements gives our proof for pseudo-metric spaces.}, and only note that it uses convergence in probability, which we had to generalize in \mathlib from pseudo-metric to extended pseudo-metric spaces.
From a distance zero modification $Y$ of $X$, we can obtain a true modification in a pseudo-metric space as follows: for each $t \in T$, let $Z_t$ be equal to $X_t$ on the event where $d_E(X_t, Y_t) = 0$, and equal to $Y_0$ otherwise.
Then $Z_t = X_t$ almost surely, and since $d_E(Z_t(\omega), Y_t(\omega)) = 0$ for all $t$ and all $\omega$, $Z$ has the same continuity properties as $Y$.
Finally, $Z$ is measurable (that is, $Z_t$ is measurable for all $t$), since a random variable that is at distance zero from a measurable random variable is itself measurable.

Now that we have a modification with Hölder continuous paths of order $\beta$ for a fixed $\beta \in (0, (q - d)/p)$, we want to obtain a modification with Hölder continuous paths of order $\gamma$ for all $\gamma \in (0, (q - d)/p)$ simultaneously.
Since $T$ has finite diameter, being Hölder continuous of order $\beta$ implies being Hölder continuous of order $\gamma$ for all $\gamma \in (0, \beta]$ (which we proved at \bmlink{Auxiliary/Topology}{L85-L87}{HolderOnWith.mono_right}).

Let $(\beta_n)_{n \in \mathbb{N}}$ be an increasing sequence in $(0, (q - d)/p)$ converging to $(q - d)/p$.
For each $n$, we obtain a modification $Y^n$ with Hölder continuous paths of order $\beta_n$.
Those $Y^n$ are all distance zero modifications of each other and are continuous, hence indistinguishable (in the sense of having distance zero).
Thus almost surely, for all $t \in T$, $d_E(Y^n_t, Y^1_t) = 0$ for all $n$.
We can thus define a process $Z$ by $Z_t(\omega) = Y^1_t(\omega)$ if $\omega$ is in the event where all $Y^n$ are at distance zero, and $Z_t(\omega) = x_0$ otherwise.
Then $Z$ is a modification of $X$ and has Hölder continuous paths of order $\beta_n$ for all $n$.
It thus has Hölder continuous paths of order $\gamma$ for all $\gamma \in (0, (q - d)/p)$.

That discussion is very close to the one in \cite{kratschmer2023kolmogorov} (up to our use of distance zero instead of equality), and we omitted an aspect that proved more subtle than expected during the formalization: the measurability of $Z$.
The process $Z$ is defined as being equal to one of two measurable processes depending on whether $\omega$ belongs to a certain event or not.
In order to ensure that $Z$ is measurable, that event must be measurable.
In particular we need the event $\{\omega \mid \forall t \in T, \: d_E(Y_t^1(\omega), Y_t^n(\omega)) = 0\}$ that two processes are at distance zero to be measurable.
We use continuity to reduce to a countable subset of $T$, and it suffices then to prove that the set $\{\omega \mid d_E(Y_t^1(\omega), Y_t^n(\omega)) = 0\}$ is measurable for all $t \in T$.
This is true if the distance $d_E(Y_t^1, Y_t^n)$ is a measurable function, which would follow from the measurability of the pairs $(Y_t^1, Y_t^n)$ with respect to $\mathcal{B}(E \times E)$.
However, as currently stated, Theorem~\ref{thm:kolmogorov_chentsov} does not guarantee that the two modifications $Y^1$ and $Y^n$ satisfy that property (unless $E$ is second-countable, in which case that follows from the measurability of each $Y^i$).
The particular construction of $Y$ we used in the proof will be important to ensure the measurability: this phenomenon of relying on the particular object constructed in the proof rather than purely on the statement of the theorem that gives existence is not ideal but also not uncommon in mathematics and is an often encountered challenge in formalization \cite{morelTalk}.
To solve it, we introduced the property for a process $Y$ of being at distance zero of a limit of $X$ on an event of probability one, and at distance zero of a fixed arbitrary value outside of that event (\bmlink{Continuity/KolmogorovChentsov}{L557-L564}{ProbabilityTheory.IsLimitOfIndicator}), and our formal statement of Theorem~\ref{thm:kolmogorov_chentsov} ensures that the modification we build satisfies that property.
The fact that all $Y^i$ satisfy that property is sufficient to ensure the measurability of $Z$ and to prove that $Z$ itself has the same property, which can then be used in further constructions (needed for the localized Theorem~\ref{thm:localized_holder_modification_sup}).

\section{Construction of a Brownian motion and Wiener measure}
\label{S:BM}

All the results presented so far can be used together to construct a Brownian motion indexed by $\mathbb{R}_+$ with values in $\mathbb{R}$.
We build a stochastic process $(B_t)_{t \in \mathbb{R}_+}$ such that
\begin{itemize}
  \item $B$ has continuous paths,
  \item $B_0 = 0$ almost surely,
  \item for all $n \in \mathbb{N}$ and $t_1, \ldots, t_n \in \mathbb{R}_+$, the vector $(B_{t_1}, \ldots, B_{t_n})$ has a multivariate normal distribution with mean $0$ and covariance matrix given by $\mathbf{cov}(B_{t_i}, B_{t_j}) = t_i \wedge t_j$.
\end{itemize}

We build the finite-dimensional distributions of $B$ using the multivariate normal distribution defined in Section~\ref{S:gaussian}: the prescribed covariance matrix is positive semi-definite since it can be written as a Gram matrix, which means that the multivariate normal distribution with mean 0 and that covariance matrix is well-defined.
We then check that they form a projective family and apply the Kolmogorov extension theorem (Theorem~\ref{T1}) to obtain a measure $P_B$ on $\mathbb{R}^{\mathbb{R}_+}$ which is the projective limit of that family.

Let $(C_t)_{t \in \mathbb{R}_+}$ be the canonical process on $\mathbb{R}^{\mathbb{R}_+}$, defined by $C_t(\omega) = \omega(t)$ for all $\omega \in \mathbb{R}^{\mathbb{R}_+}$.
The process $C$ is measurable, and its finite-dimensional distributions in the measure space $(\mathbb{R}^{\mathbb{R}_+}, P_B)$ are the projective family we defined.

Finally, we check that $C$ satisfies the Kolmogorov condition for exponents $(2n, n)$ for all $n \in \mathbb{N}$ with constant $M = (2n - 1)!!$ (the double factorial, product of all odd integers up to $2n - 1$).
We can then apply Theorem~\ref{thm:localized_holder_modification_sup} to obtain a modification $B$ of $C$ with locally Hölder continuous paths of all orders $\gamma \in (0, 1/2)$.
In particular, $B$ has continuous paths.
The process $B$ is a Brownian motion indexed by $\mathbb{R}_+$ with values in $\mathbb{R}$.

\paragraph{Properties of the Brownian motion}

Following the ideas exposed in Section~\ref{S:RV-formalism}, we do not just want to build an actual process which is a Brownian motion, we also want to be able to state that some stochastic process is a Brownian motion. Thus, as previewed in the introduction, we introduced a predicate \lean{IsBrownian X P}, similar to \lean{IsGaussianProcess X P}.

\begin{definition}[\bmlink{Gaussian/BrownianMotion}{L490-L493}{ProbabilityTheory.IsBrownian}]\label{def:IsBrownian}
A stochastic process $X$ on $\mathbb{R}_+$ with values in $\mathbb{R}$ is said to be a \emph{Brownian motion} if $X$ has continuous paths almost surely and for all $n \in \mathbb{N}$ and $t_1, \ldots, t_n \in \mathbb{R}_+$, the vector $(X_{t_1}, \ldots, X_{t_n})$ has a multivariate normal distribution with mean $0$ and covariance matrix given by $\mathbf{cov}(X_{t_i}, X_{t_j}) = t_i \wedge t_j$.
\end{definition}

The process $B$ we constructed satisfies that property, with continuous paths (not just almost surely).
To check the soundness of our formalization, we then formalized a few basic properties of the Brownian motion, which one can for example find in \cite{LeGall2016}.

Definition~\ref{def:IsBrownian} is not necessarily the most convenient way to prove that a stochastic process is a Brownian motion. Therefore we also formalized two characterizations that we present now. The first one is Proposition~\ref{prop:isBrownianOfCov}.

\begin{proposition}[\bmlink{Gaussian/BrownianMotion}{L355-L357}{ProbabilityTheory.IsGaussianProcess.isPreBrownian_of_covariance}]\label{prop:isBrownianOfCov}
A centered and almost surely continuous Gaussian process $X$ such that for all $s \le t$, $\mathbf{cov}(X_s, X_t) = s$ is a Brownian motion.
\end{proposition}

This is certainly the easiest way to prove that a given process is Gaussian. It allowed us to very easily prove the following basic properties of the Brownian motion.

\begin{proposition}
Consider $B$ a Brownian motion. Then:
\begin{enumerate}[noitemsep]
  \item for $c>0$, the process $B_{ct}/\sqrt{c}$ is a Brownian motion;
  \item \textbf{weak Markov property:} for $t_0 \in \mathbb{R}_+$, the process $B_{t_0 + t} - B_{t_0}$ is a Brownian motion independent of $(B_t)_{t \le t_0}$;
  \item $t B_{1/t}$ is a Brownian motion;
  \item $B_t/t$ tends to 0 almost surely as $t$ tends to infinity.
\end{enumerate}
\end{proposition}

Let us focus on the first three statements, the fourth one being an easy consequence of the third. Computing the expectation and the variance is quite easy. On paper, proving that these are Gaussian processes is also straightforward, because a linear combination of marginals of any of theses processes is also clearly a linear combination of marginals of the process $B$, and thus has Gaussian law. This is not as obvious to formalize, because the notion of "being a linear combination of marginals" is a bit messy to work with. What we were able to formalize nicely however is the following lemma.

\begin{lemma}[\bmlink{Gaussian/GaussianProcess}{L217-L223}{ProbabilityTheory.IsGaussianProcess.of_isGaussianProcess}]\label{lem:isGaussOfIsGauss}
Let $(X_t)_{t \in T}$ be an $E$-valued Gaussian process and $(Y_s)_{s \in S}$ be an $F$-valued stochastic process. If for all $s \in S$, there exists $I$ a finite subset of $T$ and $L : E^I \to F$ a continuous linear map such that $Y_s = L((X_t)_{t \in I})$, then $Y$ is a Gaussian process.
\end{lemma}

This statement makes proving Gaussianity in (1), (2) and (3) straightforward, especially because proving linearity is easy thanks to the different automation tactics available in \mathlib.
It is even more impressive in the case of (2). To prove independence, we rely on Lemma~\ref{lem:indepGaussProc}, which requires to prove that the joint process is Gaussian. Talking about linear combinations in the setting of the joint process would have been hard, but the phrasing of Lemma~\ref{lem:isGaussOfIsGauss} trivializes the formalization.

It remains to discuss the continuity. It is immediate in (1) and (2), and in (3) it is also obvious on $(0, +\infty)$, but the continuity at $0$ is more involved. Because we proved that the process has the law of the Brownian motion, we can apply Theorem~\ref{thm:kolmogorov_chentsov} to obtain a Brownian motion which is a modification of $t B_{1/t}$. Then we can exploit continuity to conclude. To formalize this argument, we took inspiration in the way measurable functions are formalized in \Lean.

Integration is defined for measurable functions, but if $f$ is measurable and $g$ is almost everywhere equal to $f$, then it also makes sense to integrate $g$, even though $g$ is not necessarily measurable. This motivates the notion of \emph{almost-everywhere measurability} (\mllink{MeasureTheory/Measure/MeasureSpaceDef}{L403-L407}{AEMeasurable}), which proved to be a very convenient way to manipulate almost everywhere equality in \Lean. In particular, if $f$ is an almost-everywhere measurable function, and this fact is registered in \Lean as an hypothesis \lean{hf}, then \lean{hf.mk f} is a measurable function that is almost everywhere equal to $f$. This mechanism allows to transport properties of measurable functions to almost everywhere measurable functions.

Mimicking this process, we introduced a predicate \lean{IsPreBrownian X P}, which registers that the finite dimensional laws of $X$ are those of a Brownian motion, but without assuming continuity. Then, building on Theorem~\ref{thm:kolmogorov_chentsov}, we defined \bmlink{Gaussian/BrownianMotion}{L298-L300}{IsPreBrownian.mk}, which to a pre-Brownian motion associates a Brownian motion satisfying all the properties we can expect. In particular it is \emph{everywhere} locally $\alpha$-Hölder when $\alpha \in (0, 1/2)$.

The second characterization of the Brownian motion we formalized is the one below.

\begin{proposition}[\bmlink{Gaussian/BrownianMotion}{L391-L393}{ProbabilityTheory.HasIndepIncrements.isPreBrownian_of_hasLaw}]\label{prop:isBrownianOfIndep}
An almost surely continuous stochastic process with independent increments such that for all $t$, $X_t$ has law $N(0,t)$ is a Brownian motion.
\end{proposition}

This characterization required us to define a predicate \lean{HasIndepIncrements X P} stating that the process $X$ has independent increments under the measure $P$.
This proved less straightforward than initially anticipated.
The main issue is that we manipulate finite-dimensional laws as distributions over unordered families. However, stating that a process has independent increments requires an ordering. Therefore exploiting our definition of \lean{HasIndepIncrements} to prove results about finite-dimensional distributions required us to develop a small set of lemmas to link the unordered families with the ordered ones.

In particular, the way to prove Proposition~\ref{prop:isBrownianOfIndep} is to first show that the hypotheses imply that the process is Gaussian, and then apply Proposition~\ref{prop:isBrownianOfCov}. To do that we prove this more general result:

\begin{lemma}[\bmlink{Gaussian/BrownianMotion}{L161-L168}{ProbabilityTheory.HasIndepIncrements.isGaussianProcess}]
Let $E$ be a second-countable Banach space and $T$ be an index set with a least element $\bot$. Consider $(X_t)_{t \in T}$ an $E$-valued stochastic process with independent increments, and such that for all $t$, $X_t$ has a Gaussian law. Further assume that $X_\bot$ is almost surely equal to $0$. Then $X$ is a Gaussian process.
\end{lemma}

This is probably the lemma which required most work. The idea is to say that $(X_{t_1}, ..., X_{t_n})$ can be expressed as a linear combination of the increments $(X_{t_1} - X_\bot, ..., X_{t_n} - X_{t_{n-1}})$. Because the increments are independent and each $X_t$ is Gaussian, we can prove that each increment is Gaussian. Because they are independent, their linear combination is also Gaussian, and so $(X_{t_1}, ..., X_{t_n})$ is Gaussian. Making the linear combination explicit and proving the equality was rather technical.

To conclude, let us mention a last result which is fundamental to prove more advanced properties of the Brownian motion, such as the almost sure nowhere differentiability.

\begin{proposition}[Blumenthal's zero-one law]
Let $B$ be a Brownian motion and denote by $\mathcal{F}$ the canonical filtration. The $\sigma$-algebra $\bigcap_{t > 0} \mathcal{F}_t$ is trivial, i.e.\ it contains only sets which have probability zero or one.
\end{proposition}

We followed the proof from \cite{LeGall2016}. This led to proving important results about independence, mainly that a process $(X_t)_{t \in T}$ is independent form a $\sigma$-algebra $\mathcal{A}$ if for all $t_1, ..., t_n \in T$, for all bounded continuous function $f$ and for all $A \in \mathcal{A}$, $P[\mathbb{I}_A f(X_{t_1}, ..., X_{t_n})] = P(A) * P[f(X_{t_1}, ..., X_{t_n})]$. We also relied on the weak Markov property which proof was described above.

Besides providing useful properties of the Brownian motion, those results also serve as verification of the fact that our Lean definition actually corresponds to the usual notion of Brownian motion. It also illustrates that working with our definitions in a formalization setting is practical, which is essential if we hope to formalize more advanced results.

\paragraph{Wiener measure}

The Brownian motion $B$ can be seen as a random variable on the space $\mathbb{R}^{\mathbb{R}_+}$, with the product $\sigma$-algebra and the measure $P_B$.
It also has continuous paths, and thus is a random variable on the subspace of continuous functions of $\mathbb{R}^{\mathbb{R}_+}$, with the subset $\sigma$-algebra coming from the product $\sigma$-algebra on $\mathbb{R}^{\mathbb{R}_+}$.

We may want to consider instead the space of continuous functions $\mathcal{C}(\mathbb{R}_+, \mathbb{R})$, equipped with the Borel $\sigma$-algebra generated by the topology of uniform convergence on compact sets.

We prove in \bmlink{Gaussian/BrownianMotion}{L813-L815}{ProbabilityTheory.ContinuousMap.borel_eq_iSup_comap_eval} that those two $\sigma$-algebras on continuous functions from a topological space $X$ to another $Y$ are equal whenever $X$ is a second-countable, locally compact space and $Y$ is second-countable and regular, equipped with the Borel $\sigma$-algebra.
These assumptions are satisfied for $X = \mathbb{R}_+$ and $Y = \mathbb{R}$.
Thus $B$ is a random variable on $(\mathcal{C}(\mathbb{R}_+, \mathbb{R}), \mathcal{B}(\mathcal{C}(\mathbb{R}_+, \mathbb{R})))$.
The law of $B$ in that space is called the Wiener measure (\bmlink{Gaussian/BrownianMotion}{L970-L971}{ProbabilityTheory.wienerMeasure}).

\section*{Acknowledgments}

We would like to thank Markus Himmel for his many contributions to the formalization of the Kolmogorov-Chentsov proof, as well as Sébastien Gouëzel, whose suggestions helped us greatly improve Theorem~\ref{thm:kolmogorov_chentsov} by removing any need for a metric space (instead of pseudo-metric space) assumption.
We are also grateful to Jonas Bayer, Lorenzo Luccioli, Alessio Rondelli, Jérémy Scanvic for their contributions to various parts of the formalization and to Pietro Monticone for his help with technical aspects of setting up the project.

This project would not have been possible without the foundational work of the \Lean and \mathlib community.
We would also like to thank the \mathlib reviewers, who provided very useful feedback on the parts of the formalization that we incorporated to \mathlib. PP thanks the Vector foundation for partial funding (LeanAI, P2022-0241).


\printbibliography

\end{document}